\documentclass[11pt]{article}
\usepackage[margin=1in]{geometry}

\usepackage{graphicx}%
\usepackage{tikz}
\usetikzlibrary{arrows.meta,decorations.pathreplacing}
\usepackage{multirow}%
\usepackage{amsmath,amssymb,amsfonts}%
\usepackage{xcolor}%
\usepackage{amsthm,stmaryrd,mathtools}
\usepackage{subcaption}
\usepackage{enumitem}

\usepackage{algorithm}
\usepackage{algorithmic}
\usepackage[hidelinks]{hyperref}

\graphicspath{ {./figs/}}

\newtheorem{theorem}{Theorem}
\newtheorem{proposition}[theorem]{Proposition}%
\newtheorem{lemma}[theorem]{Lemma}

\title{Operator Learning on the Data-Driven Multiscale Space \\ for Nonlinear Flow in Random Heterogeneous Porous Media}

\author{
Maria Vasilyeva 
\thanks{Department of Mathematics and Statistics, Texas A\&M University-Corpus Christi,   Corpus Christi, TX, USA\\
Email: {\tt maria.vasilyeva@tamucc.edu}.}
\and
Yiran Wang
\thanks{Department of Mathematics, The University of Alabama, Tuscaloosa, AL, USA\\
Email: {\tt ywang348@ua.edu}.}
\and
Raphael Pangilinan
\thanks{Department of Computer Science, Texas A\&M University-Corpus Christi,   Corpus Christi, TX, USA\\
Email: {\tt rpangilinan@islander.tamucc.edu}.}
}

\date{}

\begin{document}
\maketitle

\begin{abstract}
We present a neural operator surrogate for nonlinear flow in random high-contrast porous media, in which the solution is represented in a local multiscale space constructed from fine-grid snapshots. The fixed local modes define the trunk, and a branch network maps the encoded permeability to the solution coefficients and replaces the nonlinear Galerkin coarse solve. The prediction can be transferred to other grids without retraining. The error splits into a projection error and a learning error, and one correction step gives computable indicators of both. For a Forchheimer-type model with pressure-dependent permeability, the space with five local modes gives smaller mean relative $L^2$ errors than the Galerkin solution in both tested families of media, and the sparse local basis has a low storage cost. In four representative cases, a few fine-grid corrections reach the prescribed error targets faster than the reference fine-grid solve, which is computed to a higher accuracy.
\end{abstract}


\section{Introduction}
The simulation of nonlinear flow in heterogeneous porous media requires fine grids that resolve the channels, barriers and high-contrast variations of the permeability, and for random media the problem must be solved for many realizations. Two classical ways to reduce this cost are multiscale methods and model order reduction. Multiscale methods transfer the fine-scale information to a coarse grid, on which the problem is then solved. Homogenization and numerical upscaling compute effective coefficients from local problems \cite{bensoussan1978asymptotic, durlofsky1991numerical}. Multiscale methods construct coarse basis functions that contain the fine-scale variation of the solution, for example the multiscale finite element method and its generalizations \cite{hou1997multiscale, efendiev2009multiscale, efendiev2013generalized, chung2014adaptive, chung2018constraint}, the multiscale finite volume method \cite{jenny2003multi}, the variational multiscale method \cite{hughes1998variational}, the heterogeneous multiscale method \cite{e2003heterogeneous}, the multiscale mortar method \cite{arbogast2007multiscale}, and the localized orthogonal decomposition \cite{malqvist2014localization}. These methods build the coarse model from local problems for a given medium. Model order reduction takes the complementary route and extracts a global low-dimensional space from fine-scale solutions computed for many realizations, an approach that is widely used in reservoir simulation \cite{ghasemi2016model}. Among such methods, proper orthogonal decomposition (POD) identifies the dominant modes of the snapshot space and gives an optimal low-rank approximation of the solution data \cite{quarteroni2015reduced, chatterjee2000introduction, benner2017model}, and reduced basis methods project the governing equations onto these modes and solve the resulting small system for each new realization \cite{quarteroni2011certified, hesthaven2016certified}. For nonlinear problems, the nonlinear terms need an additional approximation such as the discrete empirical interpolation method \cite{chaturantabut2009discrete}, and the reduced basis can be adapted locally in the parameter space \cite{amsallem2012nonlinear, ghasemi2016localized}. In these methods a reduced nonlinear system is still assembled and solved online for every realization.

Operator learning provides an alternative that avoids the online solution of the reduced system. A neural network learns the solution operator that maps material properties, source terms, boundary conditions or other inputs to the solution field, and it is evaluated efficiently after offline training. Neural operator architectures include DeepONet \cite{lu2021learning, kovachki2024operator}, the Fourier neural operator \cite{li2020fourier}, the wavelet neural operator \cite{tripura2022wavelet} and spectral neural operators \cite{fanaskov2023spectral}, and their accuracy and cost are compared in \cite{dehoop2022cost, lu2022comprehensive}. A central choice is the space that represents the output, the trunk of the operator network. The trunk can be learned by a convolutional autoencoder \cite{oommen2022learning}, built from trainable localized finite element basis functions \cite{zhang2025finite}, or learned as continuous nonlinear basis functions by Neural-POD \cite{mou2026neural}, and mixtures of experts \cite{jacobs1991adaptive, shazeer2017outrageously} add locality through a partition-of-unity blend of local trunks \cite{sharma2024ensemble}, experts on regions of the input space \cite{kratsios2024mixtures}, or a learned domain decomposition \cite{stiller2020large}. The link to model order reduction is direct when the trunk is a fixed global POD basis and the network predicts the corresponding coefficients \cite{lu2021learning, lu2022comprehensive, bhattacharya2021model, wang2025reduced}. Such global representations can require substantial storage to capture localized features of the pressure and flux fields in high-contrast media with channels and barriers.

Local multiscale spaces are a natural representation for such media, and machine learning has been combined with multiscale methods in two main ways. The first way learns the ingredients of the multiscale discretization for each new medium. Neural networks predict the local multiscale basis functions and the local coarse-scale matrices from the permeability \cite{wang2019prediction}, the localized spectral subspaces of the multiscale space \cite{rudikov2025locally} and the multiscale prolongation operator of a two-grid preconditioner \cite{liu2024learning}. The same idea is applied to the local coarse-scale matrices of numerical homogenization, with an analysis of the resulting error \cite{kropfl2022operator, kropfl2023neural}, to the local Dirichlet-to-Neumann operators of a multiscale hybridizable discontinuous Galerkin method \cite{haines2026neural}, to coarse-grid models of the flow dynamics \cite{wang2020deep}, and to the upscaled parameters of finite volume and finite element approximations in our previous work \cite{vasilyeva2020learning, vasilyeva2021machine}. The second way learns the reduced solution itself. In \cite{wang2022deep}, a global POD space is computed from multiscale solutions for one representative permeability field, and a small network maps a few Karhunen--Lo\`eve coefficients of the permeability, and optionally the time, to the POD coefficients of the solution of a parabolic flow problem. Related modular constructions train a network on a small subdomain and assemble its predictions iteratively on larger domains \cite{wang2022mosaic}, or train data-driven elements locally and couple them as in the finite element method \cite{choi2025defining, chung2026latent}, an approach proposed as a route towards foundation models for computational science. A complementary idea is a universal space that is computed once and reused for all realizations. Such a space was obtained in \cite{vasilyeva2021preconditioning} by collecting multiscale basis functions over many realizations of the stochastic field and extracting a compact set of modes with POD.

In this work, we combine these ideas and propose \textit{operator learning on a local data-driven multiscale space} for nonlinear flow in random heterogeneous porous media. The space is constructed once from fine-grid solution snapshots by local data compression on the coarse neighborhoods and partition-of-unity coupling, so that each basis function is supported on one neighborhood and represents one of its dominant flow patterns. The local support makes the basis matrix sparse and reduces its storage cost. In the examples considered here, the local spaces are more accurate at comparable basis-storage budgets, whereas global POD can be more accurate at comparable numbers of coefficients. On this space, a neural operator replaces the Galerkin coarse solve and maps the sampled permeability field directly to the coefficients of the solution, with no reduced nonlinear solve and no approximation of nonlinear terms. The coefficient-map dimension is independent of the output fine grid, and the same coefficients can be decoded on other output grids by transferring the fixed basis functions. The output-grid transfer does not introduce new medium information or scales absent from the training-grid basis. The prediction also serves as the starting state of a few fine-grid correction steps, which recover fine-scale information that is not contained in the reduced space. 
The error of the prediction splits into a projection error and a learning error, and one correction step gives computable indicators of both. 
In the numerical experiments, the neural operator is compared with the $L^2$ projection and the Galerkin coarse solve on the same spaces and with global POD spaces.

The paper is organized as follows.
Section~\ref{sec:model} introduces the nonlinear flow model and the random heterogeneous datasets.
Section~\ref{sec:ms} constructs the local data-driven multiscale space and the corresponding Galerkin coarse solve.
Section~\ref{sec:no} develops operator learning on the multiscale space and describes the prediction on other grids and the correction of the prediction.
Section~\ref{sec:res} presents numerical results.
The conclusions follow in Section~\ref{sec:conclusions}. 
The global POD space is described in Appendix~\ref{AppA}, and the error analysis is given in Appendix~\ref{AppB}.

\section{Problem Formulation}
\label{sec:model}

We first introduce the isotropic vector Forchheimer law with pressure-dependent permeability (Section~\ref{sec:flow}). Section~\ref{sec:fine} specifies the conservative face-normal finite-volume model used to define the reported discrete problems and explains its distinction from the isotropic vector law. Section~\ref{sec:data} describes the random permeability fields and the dataset of fine-grid solutions.

\subsection{Nonlinear Forchheimer Flow}
\label{sec:flow}

We consider steady-state single-phase flow in a porous medium in $\Omega = [0,L_x] \times[0,L_y]$. 
For an incompressible fluid the flow is governed by the mass conservation equation
\[
\nabla \cdot q = f \quad \text{in } \Omega,
\]
where $q$ is the flux, and $f$ is a source or sink term.
For slow flow the flux is given by Darcy's law
\[
q = - k(x, p) \nabla p, \quad k(x, p) = \kappa(x, p)/\mu, 
\]
where $\mu$ is the fluid viscosity and $\kappa$ is the heterogeneous permeability. Here we consider a pressure-dependent permeability of the exponential form
\[
\kappa(x, p) = \kappa_0(x) \exp\big(\alpha (p-p_0)\big), 
\]
where $\kappa_0$ and $p_0$ are the reference permeability and pressure, and $\alpha$ controls the strength of the exponential nonlinearity. Such a law describes stress-sensitive media, in which a higher fluid pressure reduces the effective stress on the solid skeleton, opens the pores and microfractures and increases the permeability.

For larger velocities the inertial effects are taken into account by the Forchheimer law
\[
-\nabla p = \frac{\mu}{\kappa(x,p)} q
+ \beta \rho |q| q,
\]
where $\rho$ is the density, and $\beta$ is the Forchheimer coefficient, which depends on the permeability as $\beta=\beta(x,p)=\beta_0/\sqrt{\kappa(x,p)}$ with a constant $\beta_0$. 

Taking the norm of both sides gives a quadratic equation for $|q|$, and its positive root gives \cite{aulisa2009analysis, arraras2019geometric, vasilyeva2026implicit}
\[
G(x,p,|\nabla p|)\,q = -\nabla p, 
\quad 
G(x, p,s)
=
\frac{1}{2}
\left(
\frac{\mu}{\kappa(x, p)}
+
\sqrt{
\left(\frac{\mu}{\kappa(x, p)}\right)^2
+
4\beta\rho\,s
}
\right),
\]
where the flux is $q=-\lambda(x,p,|\nabla p|)\nabla p$ and $\lambda(x,p,s)=1/G(x,p,s)$ is the nonlinear mobility. For $\beta_0=0$ we have $\lambda=\kappa/\mu$, and the relation reduces to Darcy's law.
Next, we substitute the flux relation into the conservation law, resulting in the nonlinear
pressure equation
\begin{equation}
\label{eq:pde}
- \nabla \cdot
\left(
\lambda(x, p,|\nabla p|)\nabla p
\right)
= f
\quad \text{in } \Omega, 
\end{equation}
with the Dirichlet boundary condition $p=a$ on $\partial\Omega$, where $a$ is a constant. The problem has two sources of nonlinearity, the dependence of the permeability on the pressure and the dependence of the mobility on the pressure gradient. 

\subsection{Discretization and numerical solution}
\label{sec:fine}

We define the discrete problem by conservative cell-centred finite-volume fluxes on a Cartesian fine grid $\mathcal{T}_h$ with $N_x\times N_y$ cells of size $h_x\times h_y$, where $h_x=L_x/N_x$, $h_y=L_y/N_y$ and $N=N_xN_y$. The fine grid resolves the heterogeneous permeability field. The cell $K_I$ in the column $i$ and the row $j$ has the global index $I=(j-1)N_x+i$, $i=1,\ldots,N_x$, $j=1,\ldots,N_y$, its centre is denoted by $x_I$, and $p_I\approx p(x_I)$ is the cell pressure. Let $N(I)$ be the set of the cells that share a face with $K_I$. For $J\in N(I)$ let $E_{IJ}$ be the common face of $K_I$ and $K_J$, $|E_{IJ}|$ its length and $d_{IJ}$ the distance between $x_I$ and $x_J$. The flux through $E_{IJ}$ is approximated by $T_{IJ}(p)\,(p_I-p_J)$ with the transmissibility
\[
T_{IJ}(p)=\lambda_{IJ}\,\frac{|E_{IJ}|}{d_{IJ}},\quad
\lambda_{IJ} = 2 \left( \dfrac{\mu}{\kappa_{IJ}}+\sqrt{\Big(\dfrac{\mu}{\kappa_{IJ}}\Big)^2+4\beta_{IJ}\rho\,\dfrac{|p_I-p_J|}{d_{IJ}}} \right)^{-1} ,
\]
which is the mobility $\lambda$ evaluated with the harmonic average $\kappa_{IJ}=2\kappa_I\kappa_J/(\kappa_I+\kappa_J)$ of the cell permeabilities $\kappa_I=\kappa(x_I,p_I)$, with $\beta_{IJ}=\beta_0/\sqrt{\kappa_{IJ}}$ and with the face-normal difference quotient $|p_I-p_J|/d_{IJ}$ as the argument of the mobility. For a face $E_{ID}$ of $K_I$ on $\partial\Omega$ the flux is $T_{ID}(p)\,(p_I-a)$, where $T_{ID}$ is defined in the same way with $\kappa_I$ in place of $\kappa_{IJ}$, $p_J$ replaced by $a$ and $d_{ID}$ equal to the distance from $x_I$ to the face. With $f_I$ the value of $f$ in $K_I$ and $|K_I|=h_xh_y$, the discrete conservation law reads
\[
\sum_{J\in N(I)}T_{IJ}(p)\,(p_I-p_J)
+\sum_{E_{ID}\subset\partial\Omega}T_{ID}(p)\,(p_I-a)
=|K_I|\,f_I,\quad I=1,\ldots,N,
\]
and in matrix form
\begin{equation}
\label{eq:fine}
A(p)\,p=b(p),
\end{equation}
where $p=\{p_I\}_{I=1}^{N}\in\mathbb{R}^N$ is the vector of the cell pressures, the matrix $A(p)$ contains the transmissibilities and the vector $b(p)$ contains the source and the boundary terms. 
For $a=0$ the vector $b$ does not depend on $p$. 

To solve the nonlinear system \eqref{eq:fine}, we employ Picard iterations. Given the current iterate $p^{n}$, we freeze the nonlinear coefficient, assemble the transmissibilities at $p^n$, and obtain the next iterate from the linear system
\begin{equation}
\label{eq:picard}
A(p^n)\,p^{n+1}=b(p^n).
\end{equation}
The iteration starts from the solution of the linear Darcy problem and is repeated until the relative change of the solution is below a prescribed tolerance.
The fine-grid solution of \eqref{eq:fine} is used both as the reference solution and to generate the dataset of Section~\ref{sec:data}. All errors in this paper are therefore measured with respect to the fine-grid solution, and the multiscale space and the neural operator inherit the discretization error of the fine grid. 

\subsection{Random Heterogeneous Porous Media and Dataset}
\label{sec:data}

We introduce a parametric description of the heterogeneous reference
permeability by writing
\[
\kappa(x,p;\zeta)
= \kappa_0(x;\zeta)\exp\big(\alpha(p-p_0)\big),
\]
where $\zeta$ denotes a realization parameter for the random heterogeneous medium. 
It may contain parameters describing random channelized permeability fields, Karhunen--Lo\`eve expansion (KLE) coefficients for exponentially correlated random media, or directly sampled realizations of heterogeneous permeability fields. For the KLE-based media the reference permeability is log-normal \cite{ghanem1991stochastic,vasilyeva2021preconditioning},
\[
\ln\kappa_0(x;\zeta)=m_{\kappa}+s_{\kappa}\,Y(x;\zeta),\quad
Y(x;\zeta)=\sum_{k=1}^{M}\sqrt{\theta_k}\,\zeta_k\,\Phi_k(x),
\]
where $(\theta_k,\Phi_k)$ are the $M$ largest eigenvalues and the corresponding eigenfunctions of the exponential covariance $C(x,y)=\sigma^2\exp\big(-\big((x_1-y_1)^2/l_x^2+(x_2-y_2)^2/l_y^2\big)^{1/2}\big)$ with the correlation lengths $l_x$ and $l_y$, and $\zeta_k$ are independent standard normal random variables. The field $Y$ is normalized to zero mean and unit variance over $\Omega$ in each realization, so that $m_{\kappa}$ and $s_{\kappa}$ are the mean and the standard deviation of $\ln\kappa_0$. For the channelized media the parameter $\zeta$ contains the number, the positions and the widths of the channels.

We generate realizations $\{\zeta_j\}_{j=1}^{N_d}$ from the prescribed distribution of heterogeneous media and, for each realization, solve the fine-grid nonlinear problem \eqref{eq:fine} to obtain the solution $p(\cdot;\zeta_j)$. This gives the dataset
\begin{equation}
\label{eq:data}
\mathcal{D} 
= 
\left\{
\kappa_0(\cdot;\zeta_j),\; p(\cdot;\zeta_j)
\right\}_{j=1}^{N_d},
\end{equation}
where permeability fields $\kappa_0(\cdot;\zeta_j)$ are the input data, while the corresponding fine-grid pressure solutions $p_j = p(\cdot;\zeta_j)$ are the output snapshots.  
These snapshots are used in two ways. First, they are used to construct the data-driven multiscale space through local POD-type compression. Second, they are used to train the coarse-grid operator that maps a new heterogeneous medium to the reduced multiscale coefficients in an operator learning framework.

\section{Data-Driven Multiscale Approximation on the Coarse Grid}
\label{sec:ms}

Next, we construct a coarse-scale approximation using a local data-driven representation.  In Section~\ref{sec:msspace} we define the coarse grid with its neighborhoods and construct the multiscale space. The local modes are computed from the snapshots of the fine-grid solutions in each neighborhood and are multiplied by partition of unity functions. In Section~\ref{sec:coarse} we define the two approximations in this space that are used in the paper, the projection of the fine-grid solution and the Galerkin coarse solution. 
Estimates of their errors are given in Appendix~\ref{sec:approx}.

\subsection{Data-Driven Multiscale Space}
\label{sec:msspace}

Let $\mathcal{T}_H$ be a Cartesian coarse grid of $\Omega$ with $N^H_x\times N^H_y$ coarse cells of size $H_x=L_x/N^H_x$ and $H_y=L_y/N^H_y$, where the coarse cells are much larger than the fine cells, $H\gg h$ ($H_x\gg h_x$ and $H_y\gg h_y$). The fine grid $\mathcal{T}_h$ is a refinement of $\mathcal{T}_H$, so that each coarse cell is the union of $(N_x/N^H_x)\times(N_y/N^H_y)$ fine cells. Let $x_I$, $I=1,\ldots,N_{\omega}$, denote the nodes of the coarse grid, $N_{\omega}=(N^H_x+1)(N^H_y+1)$. The neighborhood $\omega_I$ of the node $x_I$ is the union of the coarse cells that have $x_I$ as a vertex. It consists of four coarse cells for an interior node, of two coarse cells for a node on the boundary of $\Omega$, and of one coarse cell for a corner of $\Omega$. The neighborhoods overlap and cover $\Omega$, and we denote by $n_I$ the number of fine cells in $\omega_I$.

Let $p_j \in \mathbb{R}^{N}$ denote the fine-grid vector corresponding to $p(\cdot;\zeta_j)$, $j=1,\ldots,N_s$, where $N_s$ is the number of solution snapshots used for the construction of the space.

Next, we define a local snapshot matrix
\[
S_I = 
\left[
R_I\,p_j
\right]_{j=1}^{N_s},
\]
where $R_I$ $\in\mathbb{R}^{n_I\times N}$ is the restriction operator from the global fine grid to the degrees of freedom in the neighborhood $\omega_I$.
Then the local POD modes \cite{sirovich1987turbulence,chatterjee2000introduction} $\phi_{\ell}^{\omega_I}$ are defined by
\begin{equation}
\label{eq:lpod}
S_I S_I^T \phi_{\ell}^{\omega_I}
=
(\sigma_{\ell}^{\omega_I})^2 \phi_{\ell}^{\omega_I},
\quad
\ell = 1,2,\ldots,
\end{equation}
with the singular values ordered from largest to smallest.  
Keeping the first $m$ local modes gives the best low-rank representation of the snapshot restrictions on $\omega_I$.  In this sense, local POD keeps the dominant data patterns in each coarse neighborhood  $\omega_I$ rather than forcing a single global basis to explain all local features simultaneously.

On the Cartesian coarse grid the neighborhoods have only four different shapes, namely interior neighborhoods, boundary neighborhoods in two orientations, and corner neighborhoods. We therefore divide the neighborhoods into the corresponding shape groups $g$ and collect the local snapshots of all neighborhoods of a group in one matrix
\begin{equation}
\label{eq:spod-pool}
S_g=\left[R_I\,p_j\right]_{I\in g,\ j=1,\ldots,N_s}.
\end{equation}
One set of local modes $\phi^{g}_{\ell}$ is computed for each group from \eqref{eq:lpod} with $S_I$ replaced by $S_g$, and it is used in all neighborhoods of the group, $\phi^{\omega_I}_{\ell}=\phi^{g}_{\ell}$ for $I\in g$. Each set of modes is thus computed from $|g|N_s$ local snapshots instead of $N_s$, and only as many sets of modes as there are shapes have to be computed and stored.
The local modes are then coupled into a global discrete multiscale space by a partition of unity.

Let $\chi^{\omega_I}$ be the partition-of-unity function \cite{babuska1997partition,efendiev2013generalized} of the coarse neighborhood $\omega_I$. It is the bilinear function on the coarse grid equal to one at the node $x_I$ and to zero at all other coarse nodes. Its support is $\omega_I$, and $\sum_{I}\chi^{\omega_I}=1$ in $\Omega$. Let $\tilde{\phi}^{\omega_I}_{\ell}\in\mathbb{R}^N$ denote the extension of the local mode $\phi^{\omega_I}_{\ell}$ by zero from $\omega_I$ to the whole domain.
The data-driven multiscale basis functions then define the coarse space
\begin{equation}
\label{eq:msv}
V_{ms}
=
\operatorname{span}
\left\{
\psi^{\omega_I}_{\ell} =
\chi^{\omega_I}\,\tilde{\phi}^{\omega_I}_{\ell}
:\ 
\ell = 1,\ldots,m,\ 
I = 1,\ldots,N_{\omega}
\right\},
\end{equation}
where the product is taken cell by cell, with $\chi^{\omega_I}$ evaluated at the cell centres and $\psi^{\omega_I}_{\ell}\in\mathbb{R}^N$ is a fine-grid vector as $p$. Each basis function $\psi^{\omega_I}_{\ell}$ is normalized to unit Euclidean norm on the fine grid, and the space has $N_{ms}=mN_{\omega}$ degrees of freedom.
An approximation in this space has the form
\[
p_{ms}(\zeta)
=
\sum_{I=1}^{N_{\omega}}
\sum_{\ell=1}^{m}
U_{I,\ell}(\zeta)\,\psi^{\omega_I}_{\ell}\in\mathbb{R}^N,
\]
where the coefficients $U_{I,\ell}(\zeta)$ encode the dependence on the input
medium. 
Compared with global POD space (see Appendix~\ref{AppA}), the local data-driven multiscale space is often more efficient for random heterogeneous media because the same number of modes per patch can adapt to different local flow patterns across the domain.  It also leads to a sparse coarse representation and is well-suited for the operator-learning framework, where the neural network does not need to predict the full fine-grid pressure field, but only the reduced coefficient vector associated with the fixed local multiscale basis.

\subsection{Coarse-scale system and approximation properties}
\label{sec:coarse}

The multiscale approximation restricts the pressure to the data-driven space $V_{ms}$.  With $m$ local modes per coarse neighborhood, the number of coarse unknowns is $N_{ms}=mN_{\omega}$. 
Let $P$ be the prolongation matrix defined by local data-driven multiscale basis functions 
\begin{equation}
\label{eq:proj}
P=
\left[
\psi^{\omega_1}_1,\ldots,\psi^{\omega_1}_{m},
\ldots,
\psi^{\omega_I}_1,\ldots,\psi^{\omega_I}_{m},
\ldots,
\psi^{\omega_{N_{\omega}}}_1,\ldots,
\psi^{\omega_{N_{\omega}}}_{m}
\right], \quad  R=P^T.
\end{equation}
The basis vectors are computed once, on the fine grid $\mathcal{T}_h$ of the snapshots, and they are normalized on this grid.

A function of $V_{ms}$ is represented as $PU$ with a coefficient vector $U\in\mathbb{R}^{N_{ms}}$. For a given medium (the dependence on $\zeta$ is omitted below) we consider two approximations of the fine-grid solution $p$ in this space. The first is the projection of $p$, which requires the fine-grid solution and gives the smallest error that can be attained in the space. The second is the Galerkin coarse solution, which is obtained from the fine-grid system projected onto the space and does not require $p$.

\textit{Projection.} We assume that the columns of $P$ are linearly independent, so that $M=P^TP$ is positive definite and the projection coefficients are unique. All inverse expressions below denote linear solves. Let $\|\cdot\|$ be the Euclidean norm of fine-grid vectors. The best approximation of $p$ in $V_{ms}$ in this norm is
\begin{equation}
\label{eq:l2proj}
\Pi p=PU_*,\quad U_*=\arg\min_{V\in\mathbb{R}^{N_{ms}}}\|p-PV\|=M^{-1}P^Tp,\quad \Pi=PM^{-1}P^T,
\end{equation}
where $\Pi$ is the $L^2$-orthogonal projector onto $V_{ms}$. The projection error $\|(I-\Pi)p\|$ depends only on the space, and $\|(I-\Pi)p\|\le\|p-v\|$ for every $v\in V_{ms}$. The coefficients $U_*$ are the training targets of the neural operator in Section~\ref{sec:no}.

\textit{Galerkin coarse solve.} The coarse nonlinear problem is obtained by the Galerkin projection of the fine-grid system onto the multiscale space. The coefficients $U^H$ are defined by the condition that the residual of the fine-grid system is orthogonal to $V_{ms}$,
\begin{equation}
\label{eq:galerkin}
P^T\big(b-A(PU^H)\,PU^H\big)=0 .
\end{equation}
The same Picard linearization used on the fine grid is applied to the coarse problem
\begin{equation}
\label{eq:coarse}
A_{ms}(U^{H,n};\zeta)U^{H,n+1}=b_{ms}(\zeta),
\end{equation}
with 
\[
A_{ms}(U^{H,n};\zeta)=R A(P U^{H,n};\zeta)P,
\quad
b_{ms}(\zeta)=R b(\zeta).
\]
In each iteration the fine-grid matrix is assembled at the current approximation $PU^{H,n}$ and projected, and a linear system with $N_{ms}$ unknowns is solved. The converged Galerkin solution is denoted by $p_G=PU^H$. 
We note that the Galerkin projection framework is not restricted to the proposed data-driven local multiscale space. Other reduced spaces can also be employed to construct the prolongation and restriction operators. 
In particular, one may employ a classical global POD basis with support on the entire domain $\Omega$, or alternatively local multiscale bases whose functions are supported on coarse neighborhoods.

\section{Operator Learning on the Multiscale Space}
\label{sec:no}

In this section the repeated solution of the nonlinear coarse problem \eqref{eq:coarse} is replaced by a neural operator, which predicts the coefficients in the fixed multiscale space directly from the permeability field. Section~\ref{sec:nn} defines the neural operator with its input and output encodings, the loss function and the network. 
Section~\ref{sec:corr} corrects the prediction on the fine grid by a few steps of a global or a local correction. 
Appendix~\ref{sec:erran} gives an a priori decomposition of the error of the prediction into a projection error and a learning error, and an a posteriori estimate.

\subsection{Neural operator}
\label{sec:nn}

Classical neural operators differ in the representation of the input and of the output. DeepONet \cite{lu2021learning} represents the output in a trunk basis that is learned together with the branch network, and POD-DeepONet \cite{lu2022comprehensive} replaces the trunk by a global POD basis of the solutions. PCA-Net \cite{bhattacharya2021model} encodes both the input and the output by global principal components and learns the map between the two sets of coefficients. The Fourier neural operator \cite{li2020fourier} works with the full fields on the grid and parameterizes its layers in the Fourier space. The neural operator of this section has the encoder and decoder structure of PCA-Net and POD-DeepONet, but both bases are local and are built on the coarse neighborhoods of the multiscale method (Figure~\ref{fig:arch}).

\begin{figure}[htbp]
\centering
\begin{tikzpicture}[font=\footnotesize, >={Stealth[length=2mm]},
  fixed/.style={draw, rounded corners=2pt, fill=black!10, align=center, minimum width=1.9cm, minimum height=1.25cm},
  trained/.style={draw=blue!70!black, thick, rounded corners=2pt, fill=blue!12, align=center, minimum width=2.0cm, minimum height=1.25cm}]
\def\archeta{eta.png}\def\archpred{pred.png}\def\archref{ref.png}
\foreach \xc/\yc/\im in {1.1/0/\archeta, 13.6/0/\archpred, 13.6/-3.5/\archref} {
  \node[inner sep=0pt] at (\xc,\yc) {\includegraphics[width=2.2cm,height=2.2cm]{\im}};
  \begin{scope}[shift={(\xc-1.1,\yc-1.1)}]
    \draw[white, opacity=0.55, line width=0.15pt, step=0.275] (0,0) grid (2.2,2.2);
    \draw[black, line width=0.3pt] (0,0) rectangle (2.2,2.2);
    \draw[white, line width=1.5pt] (0.55,0.825) rectangle (1.1,1.375);
    \draw[black, line width=0.7pt] (0.55,0.825) rectangle (1.1,1.375);
  \end{scope}}
\node[align=center] at (1.1,-1.55) {log-permeability\\$\log\kappa_0(\zeta)\in\mathbb{R}^{N_{\kappa}}$};
\node[align=center] at (13.6,-1.55) {predicted pressure\\$\widehat p=P\widehat U\in\mathbb{R}^{N}$};
\node[align=center] at (13.6,-5.05) {fine-grid solution\\$p\in\mathbb{R}^{N}$};
\draw[line width=0.3pt] (0.62,0.3) -- (0.4,1.28); \node at (0.36,1.45) {$\omega_I$};
\foreach \xc/\yc/\col in {5.45/0/black!10, 9.15/0/blue!12, 9.15/-3.5/black!10} {
  \draw[fill=\col] (\xc-0.15,\yc-0.7) rectangle (\xc+0.15,\yc+0.7);
  \foreach \k in {1,...,6} \draw (\xc-0.15,\yc-0.7+0.2*\k) -- ++(0.3,0);}
\node at (5.45,1.0) {$x\in\mathbb{R}^{qN_{\omega}}$};
\node at (9.15,1.0) {$\widehat U\in\mathbb{R}^{N_{ms}}$};
\node at (9.15,-4.5) {$U_*\in\mathbb{R}^{N_{ms}}$};
\node[fixed] (E) at (3.7,0) {local PCA\\encoder $E$};
\node[trained] (G) at (7.3,0) {network\\$\mathcal{G}_{\theta}$};
\node[fixed] (P) at (11.0,0) {multiscale\\basis $P$};
\node[fixed] (Q) at (11.0,-3.5) {projection\\$(P^TP)^{-1}P^T$};
\draw[->] (2.2,0) -- (E.west); \draw[->] (E.east) -- (5.3,0); \draw[->] (5.6,0) -- (G.west);
\draw[->] (G.east) -- (9.0,0); \draw[->] (9.3,0) -- (P.west); \draw[->] (P.east) -- (12.5,0);
\draw[->] (12.5,-3.5) -- (Q.east); \draw[->] (Q.west) -- (9.3,-3.5);
\node[draw, rounded corners=7pt, inner sep=3pt, align=center] (L) at (9.15,-1.6) {loss $\mathcal{L}(\theta)$};
\draw[->] (9.15,-0.7) -- (L.north); \draw[->] (9.15,-2.8) -- (L.south);
\draw[->, dashed, blue!70!black] (L.west) -| (G.south);
\node[blue!70!black] at (7.3,-1.85) {update of $\theta$};
\draw[decorate, decoration={brace, amplitude=4pt}] (0,1.75) -- (5.6,1.75) node[midway, above=4pt] {input encoding};
\draw[decorate, decoration={brace, amplitude=4pt}] (9.0,1.75) -- (14.7,1.75) node[midway, above=4pt] {output decoding};
\node[fixed, minimum width=0.5cm, minimum height=0.3cm] at (0.3,-2.9) {};
\node[anchor=west] at (0.65,-2.9) {fixed, computed from the data before the training};
\node[trained, minimum width=0.5cm, minimum height=0.3cm] at (0.3,-3.5) {};
\node[anchor=west] at (0.65,-3.5) {trained};
\draw[->, dashed, blue!70!black] (0.05,-4.1) -- (0.55,-4.1);
\node[anchor=west] at (0.65,-4.1) {update of the parameters in the training};
\node at (10.6,-2.6) {training target};
\end{tikzpicture}
\caption{Operator learning on the multiscale space. The log-permeability is encoded by its local principal components on the coarse neighborhoods $\omega_I$, the network maps the input features $x$ to the coefficients $\widehat U$, and the pressure is reconstructed with the multiscale basis. The two encodings are fixed, and only the network is trained. The targets of the training are the coefficients $U_*$ of the projection of the fine-grid solution $p$ onto the multiscale space. 
}
\label{fig:arch}
\end{figure}

The log-permeability and the pressure are represented in these local data-driven bases, which are computed before the training and are then kept fixed, and the network is trained only to map the coefficients of the input to the coefficients of the output,
\begin{equation}
\label{eq:no}
x=E\log\kappa_0(\zeta),\quad \widehat U=\mathcal{G}_{\theta}(x),\quad \widehat p=P\,\widehat U ,
\end{equation}
where $\kappa_0(\zeta)\in\mathbb{R}^{N_{\kappa}}$ is the vector of the values of the reference permeability $\kappa_0(\cdot\,;\zeta)$ of the realization $\zeta$ (Section~\ref{sec:data}) in the $N_{\kappa}$ cells of the grid on which the permeability is given, the logarithm is taken cell by cell, $E$ is the input encoder, $\mathcal{G}_{\theta}$ is the network with the parameters $\theta$, and $P$ is the matrix \eqref{eq:proj} of the multiscale basis. The operator can be viewed as a DeepONet-like model in which $P$ is a fixed trunk and the encoder together with the network forms the branch. We describe the three parts in the order in which they are constructed.

\textit{Output encoding.} The output of the network is the vector of the $N_{ms}$ coefficients in the multiscale space $V_{ms}$ of Section~\ref{sec:ms}, and the fine-grid pressure is reconstructed by the fixed matrix $P$. For the realization $\zeta_j$ of the dataset \eqref{eq:data}, with the features $x_j=E\log\kappa_0(\zeta_j)$, the fine-grid solution $p_j=p(\cdot\,;\zeta_j)$ is encoded by its $L^2$ projection \eqref{eq:l2proj}, $U_{*,j}=M^{-1}P^Tp_j$ with $M=P^TP$, and these coefficients are the targets of the training. The network is trained in coefficient space. The objective combines the coefficient error, measured in the norm $\|V\|_M=(V^TMV)^{1/2}$, with the fixed projection residual, and normalises each sample by the squared norm of its fine-grid solution:
\begin{equation}
\label{eq:loss}
{
\mathcal{L}(\theta)
=\frac{1}{N_{\mathrm{tr}}}\sum_j\frac{(\mathcal{G}_{\theta}(x_j)-U_{*,j})^TM\,(\mathcal{G}_{\theta}(x_j)-U_{*,j})+r_j^2}{\|p_j\|^2}
=\frac{1}{N_{\mathrm{tr}}}\sum_j\frac{\|P\mathcal{G}_{\theta}(x_j)-p_j\|^2}{\|p_j\|^2},
}
\end{equation}
where $r_j=\|p_j-PU_{*,j}\|$ is the projection residual. The loss contains only vectors of size $N_{ms}$, the fixed matrix $M$ and the numbers $r_j$ and $\|p_j\|$, which are computed once before the training, so that no fine-grid field is formed during the training. The second equality follows from $P^T(p_j-PU_{*,j})=0$. It shows that the loss on the coarse grid is equal to the mean relative squared $L^2$ error of the reconstructed field on the fine grid. 
Training a network on the projection coefficients in a fixed reduced basis is the approach of the non-intrusive reduced order models \cite{hesthaven2018non} and of PCA-Net \cite{bhattacharya2021model}, where the basis is orthonormal and $M=I$. For the local multiscale basis, which is not orthogonal, the matrix $M$ is needed. The two terms are the two components of the error of the prediction. For each sample, the two squared terms in the numerator are the learning and projection errors; the projection residual is fixed by the space and has no effect on the training gradient. Appendix~\ref{sec:erran} analyses the unnormalised mean-square error and explicitly distinguishes it from the samplewise-normalised training loss \eqref{eq:loss}.

\textit{Input encoding.} The input is compressed on the same coarse neighborhoods as the output. On each neighborhood $\omega_I$ the log-permeability is restricted to the cells of the neighborhood, centred, and projected onto the first $q$ principal components,
\[
x_I=\Theta_g^T\big(R^{\kappa}_I\log\kappa_0(\zeta)-m_g\big)\in\mathbb{R}^q,\quad I\in g,
\]
where $R^{\kappa}_I$ is the restriction to the cells of the permeability grid in $\omega_I$, and the mean $m_g$ and the matrix $\Theta_g$ of the $q$ principal components are computed from the training inputs restricted to the neighborhoods of the shape group $g$. One set of components is thus computed for each neighborhood shape, as for the local modes in Section~\ref{sec:msspace}. The feature vector $x=(x_1,\ldots,x_{N_{\omega}})\in\mathbb{R}^{qN_{\omega}}$ replaces the values of the permeability in all cells. In contrast to the output, no partition of unity is used for the input. The input is only encoded and is never reconstructed from its features, so that the overlapping local representations do not have to be combined into one global function.

\textit{Network.} The map between the two sets of coefficients is a residual multilayer perceptron with $d$ layers of width $w$ and the activation $\sigma$,
\[
h_0=\sigma(W_0x+b_0),\quad h_{l}=h_{l-1}+\sigma(W_lh_{l-1}+b_l),\quad l=1,\ldots,d-1,
\]
\[
\mathcal{G}_{\theta}(x)=D_{s}(W_{d}h_{d-1}+b_{d})+m_{s},
\]
where $m_s$ and the diagonal matrix $D_s$ are the mean and the standard deviation of the target coefficients over the training set. The parameters $\theta$ are updated by a gradient-based optimizer. During the training the validation error is evaluated regularly, and the network with the smallest validation error is used.

\textit{Advantages of the two encodings.} The representation \eqref{eq:no} has the following advantages. The network maps $qN_{\omega}$ features to $N_{ms}$ coefficients. Its size and the cost of one training step do not depend on the number $N$ of fine cells. Every prediction is an element of $V_{ms}$. It can be compared with the Galerkin solution in the same space, and the same space is the coarse space of the correction in Section~\ref{sec:corr}. Since the network is trained towards the $L^2$ projection, its prediction can be more accurate in the $L^2$ norm than the Galerkin solution, which is not optimal in this norm (Section~\ref{sec:coarse}). For a fixed space, the projection error is independent of network training. The remaining learning error depends on the input encoding, network and training data, and the space itself is also constructed from data. For a new medium the prediction requires the input features, one evaluation of the network and one product with $P$. No nonlinear iteration and no assembly of the fine-grid system are needed. Both encodings are built on the coarse neighborhoods. They are computed from small local problems, and the matrix $P$ is sparse.

\textit{Prediction on other grids.} The network is trained on the grid $\mathcal{T}_h$ of the snapshots. To evaluate the prediction on another grid $\mathcal{T}_{h'}$, the basis vectors are transferred to $\mathcal{T}_{h'}$, by cell averages if $\mathcal{T}_{h'}$ is coarser and by cell-centred bilinear interpolation. 
With the resulting matrix $P_{h'}$ the field on $\mathcal{T}_{h'}$ is $\widehat p=P_{h'}\widehat U$ with the same predicted coefficients $\widehat U$, and nothing is trained again when the grid is changed. The transferred space contains only the information resolved by the basis on $\mathcal{T}_h$. This is output-grid transfer for a fixed input encoding, not a claim of generalisation to arbitrary input resolutions or newly introduced fine-scale medium features.

\subsection{Correction of the Neural Prediction}
\label{sec:corr}

The neural operator is trained with an $L^2$-based loss and gives accurate predictions in the $L^2$ norm. The energy error is not controlled by the loss, and the prediction does not satisfy the fine-grid equations exactly, especially in highly heterogeneous media. In addition, the prediction is an element of the fixed space $V_{ms}$ and cannot be more accurate than the projection of the solution onto this space. To reduce the remaining error, we combine the prediction with a small number of correction steps on the fine grid (Figure~\ref{fig:corr}). The network supplies a fast global approximation $p^0=\widehat p$, and each step acts on the residual $r^n=b-A(p^n)\,p^n$ that is left by the current approximation $p^n$.

\begin{figure}[htbp]
\centering
\begin{tikzpicture}[font=\footnotesize, >={Stealth[length=2mm]},
  box/.style={draw, rounded corners=2pt, fill=black!10, align=center, minimum height=1.15cm, inner xsep=3pt},
  trained/.style={draw=blue!70!black, thick, rounded corners=2pt, fill=blue!12, align=center, minimum height=1.15cm, inner xsep=3pt}]
\node[trained, minimum width=2.3cm] (N) at (0.5,0) {neural prediction\\$p^0=\widehat p$};
\node[box, minimum width=2.7cm] (R) at (3.7,0) {residual\\$r^n=b-A(p^n)\,p^n$\\$A_n=A(p^n)$};
\node[box, text width=5.0cm, minimum height=1.35cm] (G) at (8.4,1.05) {\textbf{global correction} \eqref{eq:lincorr}\\$A_n\,\delta p^n=r^n$ by CG with\\the AMG preconditioner};
\node[box, text width=5.0cm, minimum height=1.35cm] (L) at (8.4,-1.05) {\textbf{local correction} \eqref{eq:localcorr}\\ local problems on $K_i^{+}$
\\then a coarse step in $V_{ms}$};
\node at (8.4,0) {or};
\node[box, minimum width=2.9cm] (U) at (13.2,0) {update\\$p^{n+1}=p^n+\delta p^n$};
\draw[->] (N.east) -- (R.west);
\draw[->] (R.east) -- ++(0.35,0) |- (G.west);
\draw[->] (R.east) -- ++(0.35,0) |- (L.west);
\draw[->] (G.east) -- ++(0.35,0) |- (U.west);
\draw[->] (L.east) -- ++(0.35,0) |- (U.west);
\draw[->] (U.south) -- (13.2,-2.4) -- (3.7,-2.4) -- (R.south);
\node[above] at (8.4,-2.4) {next step, $n\to n+1$};
\end{tikzpicture}
\caption{Correction of the neural prediction on the fine grid. In every step the coefficient is frozen at the current approximation, and the correction $\delta p^n$ is computed from one linear problem, either globally or by local problems.}
\label{fig:corr}
\end{figure}

\textit{Global correction.} The residual left by the neural prediction is removed by a few steps of a global iterative solver. Starting from the neural prediction $p^0$, step $n$ freezes the coefficient at the current state, $A_n=A(p^n)$, and solves the linearised system
\begin{equation}
\label{eq:lincorr}
A_n\,\delta p^n = b-A_n p^n, \quad p^{n+1}=p^n+\delta p^n,
\end{equation}
only approximately, by the preconditioned conjugate gradient method with a loose tolerance $\varepsilon$ on the relative residual. With exact solves, \eqref{eq:lincorr} is the Picard iteration of the fine-grid solver started from the neural prediction instead of the Darcy solution. Started from the Galerkin solution $p_G$, one exact step gives $A_G^{-1}b$. Its error $z=p-A_G^{-1}b$ is the vector in the proof of Proposition~\ref{prop:galerkin}, bounded in the $A_G$ norm by $\tau\|p\|_{A_G}$. The updated solution is no longer restricted to $V_{ms}$, but it still depends on this space through the initial state $p_G$ and the frozen matrix $A_G=A(p_G)$. 
As the preconditioner we use one V-cycle of algebraic multigrid. In this setting the neural operator takes the place of the nonlinear coarse solve. At the cost of one network evaluation it supplies the state at which the coefficient is frozen, and the fine-grid work is reduced to a few linear solves of low accuracy. One step requires one assembly of $A_n$, the setup of the preconditioner and a few CG iterations, and all of them are fine-grid operations with a cost proportional to the number of fine cells.

\textit{Local correction.} The step \eqref{eq:lincorr} can also be solved approximately by local problems. Let $K_i^{+}$, $i=1,\ldots,N_K$, be the coarse cells extended by one layer of coarse cells, and let $R_i^{+}$ be the restriction of a fine-grid vector to the fine cells in $K_i^{+}$. The correction is
\begin{equation}
\label{eq:localcorr}
\delta p^n=\sum_{i=1}^{N_K}(R_i^{+})^TW_i\big(R_i^{+}A_n(R_i^{+})^T\big)^{-1}R_i^{+}\,(b-A_np^n),
\end{equation}
where $W_i$ is the diagonal matrix of the weights that account for the overlap, so that each fine cell has the weight one divided by the number of subdomains that contain it. It is followed by a coarse correction in the multiscale space, $P(P^TA_nP)^{-1}P^T$ applied to the new residual. The matrix $A_n=A(p^n)$ is assembled once in every step with the coefficients frozen at the current approximation $p^n$. All local problems are therefore linear, and the local problems are independent and can be solved in parallel.

\section{Numerical Results}
\label{sec:res}

In this section we test the method on two families of random media. Section~\ref{sec:setup} describes the data, the reduced spaces and the error measures. Section~\ref{sec:nolearn} reports the accuracy of the neural operator on the local and global reduced spaces, the storage of the spaces and the offline and online cost, compares the prediction with the projection and the Galerkin coarse solution, and gives the prediction on other grids and the computational time. Section~\ref{sec:rescorr} combines the prediction with the corrections of Section~\ref{sec:corr}.

\subsection{Model Problem and Datasets}
\label{sec:setup}

\textit{Model problem.} The domain is the unit square with homogeneous Dirichlet conditions, the source $f=10$, the viscosity $\mu=1$ and the density $\rho=10^3$. All experiments use the discrete face-normal Forchheimer-type model of Section~\ref{sec:fine}, combined with the pressure-dependent exponential permeability, with $\beta_0=10^{-4}$ and $p_0=0$, and $\alpha=10^2$ in Test~1 and $\alpha=10$ in Test~2. The dataset solutions are computed on the $240\times240$ fine grid by the Picard iteration, which starts from the Darcy solution and stops when the relative change is below $10^{-6}$, with a sparse direct solver in every iteration. The comparisons across output grids below instead use the restricted $960\times960$ solution as their common reference.

\begin{figure}[htbp]
\centering
\includegraphics[width=0.9\linewidth]{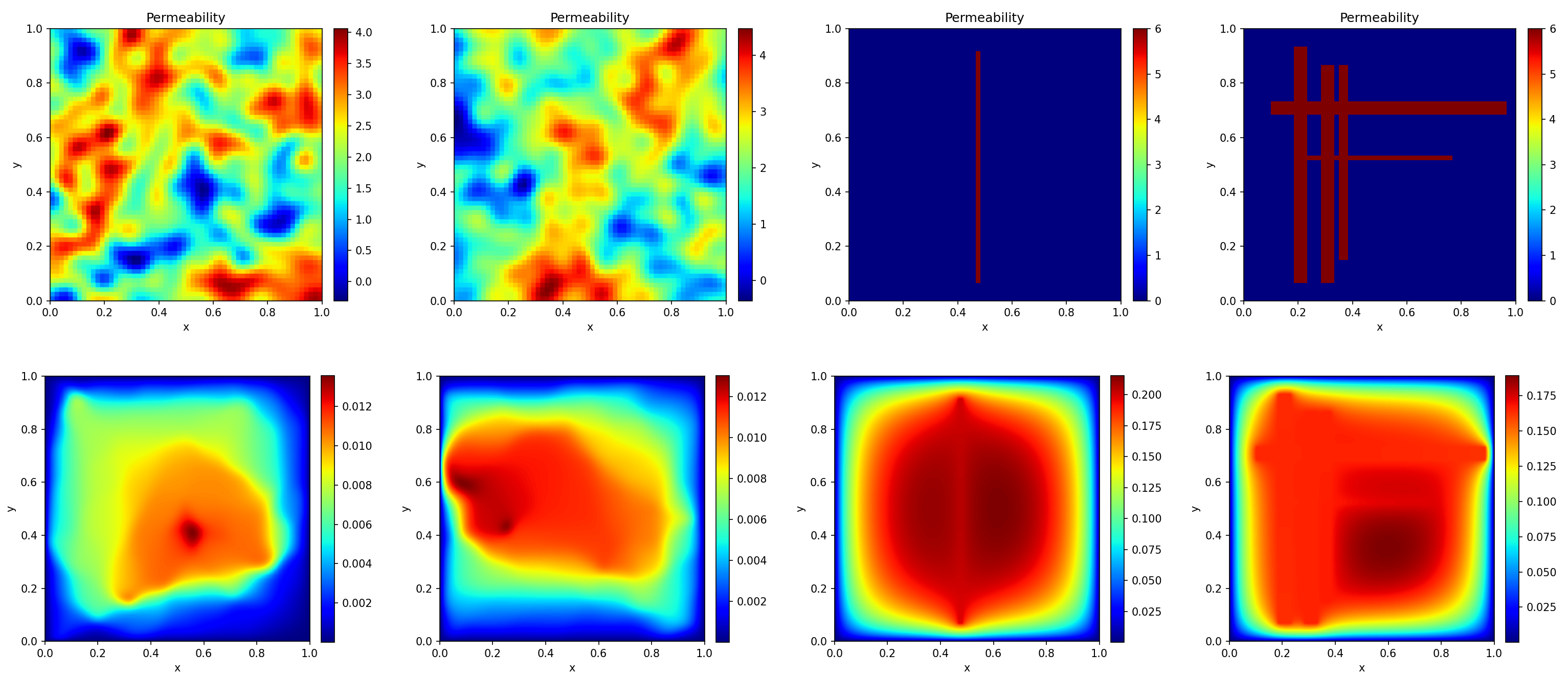}
\caption{Permeability fields and fine-grid solutions of the four separate tests. The columns show Tests 1a and 1b (KLE-based heterogeneous media) and Tests 2a and 2b (channelized media). The first row: permeability $\kappa_0$ on a logarithmic scale. The second row: fine-grid solution.}
\label{fig:sol}
\end{figure}

\textit{Permeability fields.} In Test~1 the permeability is log-normal, generated from a 300-mode Karhunen--Lo\`eve expansion with exponential covariance ($l_x=l_y=0.2$, $\sigma^2=1$, $92\%$ of the variance, $m_{\kappa}=5$, $s_{\kappa}=2$, Section~\ref{sec:data}). In Test~2 the permeability is $1$ in the background and $10^6$ in up to five horizontal and five vertical channels of random position and width. In both tests it is piecewise constant on a $60\times60$ grid. Figure~\ref{fig:sol} shows two realizations of each test (Tests 1a, 1b, 2a and 2b), which are not part of the datasets and are used as separate tests.

\textit{Datasets.} The datasets contain $17000$ random permeability fields with their fine-grid solutions for Test~1 and $12000$ for Test~2. They are split with a fixed seed into $16200$ training, $300$ validation and $500$ test realizations in Test~1, and into $11200$ training, $300$ validation and $500$ test realizations in Test~2. 
The training set is further enlarged by the eight symmetries $g$ of the square, which map the solution $p$ for the permeability $\kappa_0$ to the solution $g\cdot p$ for $g\cdot\kappa_0$. In every optimization step one symmetry is chosen at random. The validation set is used only to select the network and the test set only to report the errors, and neither is augmented.

\textit{Reduced spaces.} The coarse grid is $8\times8$, with $N_{\omega}=81$ neighborhoods in four shape groups. The local modes of each group are computed as in Section~\ref{sec:ms} from the pooled snapshots of the $4200$ training solutions. 
The space $V_{ms,1}$ has the constant function as its only local mode. 
It is a separate baseline, whereas $V_{ms,5}$, $V_{ms,10}$ and $V_{ms,20}$ are the nested local POD spaces. The global POD spaces $V_{pod,r}$ of Appendix~\ref{AppA}, $r=50,100,200,400$, are computed from the same training solutions. All spaces are fixed, and the training targets are the $L^2$ projections of the training solutions, with the mean included for the affine global POD representation.

\textit{Error measures.} The reported errors are the relative $L^2$ and energy errors in percent,
\[
\frac{\|p-p_{ms}\|}{\|p\|}\times 100\%,\quad \frac{\|p-p_{ms}\|_A}{\|p\|_A}\times 100\%,\quad \|v\|^2=v^Tv,\quad \|v\|_A^2=v^TAv,
\]
where $p$ is the fine-grid reference solution, $p_{ms}$ is the approximation, that is, the $L^2$ projection, the Galerkin solution $p_G$, the neural prediction $\widehat p$ or its corrected value, and $A=A(p_{ms})$ is the fine-grid matrix at that approximation. In the tables and figures the first three approximations in a reduced space $V$ are denoted by P+$V$, S+$V$ and NO+$V$. The $L^2$ error measures the pressure field itself. The reported energy error is a state-dependent auxiliary measure sensitive to gradients and high-contrast interfaces. For the separate tests the errors of the single realization are reported, and for the test sets the mean over the $500$ test realizations.

\subsection{Operator Learning on the Data-Driven Multiscale Space}
\label{sec:nolearn}

\textit{Networks and training.} We first evaluate the reduced spaces as the trunk of the neural operator. The networks are trained as in Section~\ref{sec:no}, with the loss in coefficient space, $q=10$ local input components per neighborhood ($810$ features). The branch network is a residual multilayer perceptron of width $1024$ with $6$ layers (GELU activation), trained by AdamW with the learning rate $5\cdot10^{-4}$, a one-cycle schedule, mini-batches of $1024$ samples and $40000$ steps, and the network with the smallest validation error is kept. The architecture and the training are the same for all spaces. The implementation is written in Python with \textit{scipy.sparse} \cite{virtanen2020scipy} and \textit{PyTorch} \cite{paszke2019pytorch}, and the computations were performed on a computer with an Apple M2 Max processor and 32 GB of memory.

\begin{table}[h!]
\centering
\footnotesize
\caption{Operator learning on global and local reduced spaces. $L^2$ / energy errors in percent. The errors of the prediction on $500$ training and on the $500$ test realizations (mean, smallest and the largest $L^2$ error), and the mean errors of the $L^2$ projection of the test solutions (Projection).}
\setlength{\tabcolsep}{2pt}
\begin{tabular}{lccccccc}
\hline
 & \multicolumn{3}{c}{Train} & \multicolumn{3}{c}{Test} & \\
 & mean & min & max & mean & min & max & Projection \\
\hline
\multicolumn{8}{c}{Test 1} \\
\hline
NO+$V_{pod,50}$ & 4.23 / 32.83 & 2.03 / 17.89 & 10.57 / 45.89 & 4.25 / 32.49 & 1.95 / 17.28 & 12.09 / 48.75 & 4.12 / 32.57 \\
NO+$V_{pod,100}$ & 2.53 / 25.10 & 1.16 / 12.34 & 7.48 / 41.58 & 2.67 / 25.15 & 1.37 / 13.97 & 7.45 / 35.90 & 2.45 / 25.19 \\
NO+$V_{pod,200}$ & 1.40 / 17.07 & 0.68 / 8.56 & 3.79 / 28.16 & 1.65 / 17.52 & 0.79 / 10.68 & 6.53 / 24.48 & 1.21 / 17.27 \\
NO+$V_{pod,400}$ & 1.02 / 12.27 & 0.54 / 7.40 & 2.44 / 18.93 & 1.31 / 12.70 & 0.54 / 7.02 & 6.47 / 17.44 & 0.47 / 9.95 \\
\hline
NO+$V_{ms,1}$ & 5.77 / 69.17 & 3.42 / 43.52 & 12.03 / 69.83 & 5.77 / 68.83 & 3.56 / 57.38 & 13.12 / 69.10 & 5.67 / 69.00 \\
NO+$V_{ms,5}$ & 1.45 / 23.36 & 0.72 / 13.08 & 6.29 / 44.00 & 1.57 / 22.91 & 0.71 / 13.41 & 6.55 / 27.93 & 1.03 / 21.76 \\
NO+$V_{ms,10}$ & 1.09 / 15.28 & 0.57 / 8.34 & 3.53 / 21.21 & 1.27 / 15.17 & 0.61 / 8.12 & 5.72 / 20.65 & 0.34 / 11.86 \\
NO+$V_{ms,20}$ & 1.02 / 12.25 & 0.53 / 6.95 & 3.43 / 16.60 & 1.20 / 12.28 & 0.60 / 8.24 & 4.99 / 16.87 & 0.12 / 6.42 \\
\hline
\multicolumn{8}{c}{Test 2} \\
\hline
NO+$V_{pod,50}$ & 1.29 / 63.66 & 0.09 / 0.82 & 2.41 / 111.63 & 1.41 / 63.10 & 0.09 / 0.82 & 4.37 / 79.89 & 1.28 / 63.14 \\
NO+$V_{pod,100}$ & 0.80 / 56.62 & 0.05 / 0.55 & 1.55 / 104.26 & 0.97 / 56.79 & 0.05 / 0.55 & 4.02 / 67.98 & 0.79 / 56.59 \\
NO+$V_{pod,200}$ & 0.52 / 49.36 & 0.04 / 0.44 & 0.87 / 84.42 & 0.72 / 50.70 & 0.04 / 0.44 & 3.92 / 60.82 & 0.47 / 49.58 \\
NO+$V_{pod,400}$ & 0.41 / 43.76 & 0.04 / 0.36 & 0.68 / 61.62 & 0.63 / 45.87 & 0.04 / 0.36 & 3.49 / 56.28 & 0.27 / 40.79 \\
\hline
NO+$V_{ms,1}$ & 5.65 / 131.98 & 4.31 / 88.82 & 7.68 / 176.73 & 5.69 / 131.98 & 4.31 / 88.82 & 7.60 / 191.32 & 5.66 / 132.00 \\
NO+$V_{ms,5}$ & 0.69 / 62.94 & 0.07 / 2.11 & 1.25 / 116.38 & 0.80 / 62.16 & 0.07 / 2.11 & 3.51 / 82.71 & 0.59 / 62.78 \\
NO+$V_{ms,10}$ & 0.48 / 53.21 & 0.03 / 0.76 & 0.78 / 76.66 & 0.63 / 52.54 & 0.03 / 0.76 & 3.96 / 66.55 & 0.31 / 49.29 \\
NO+$V_{ms,20}$ & 0.43 / 49.23 & 0.02 / 0.87 & 0.71 / 79.21 & 0.58 / 48.68 & 0.02 / 0.87 & 3.46 / 59.52 & 0.17 / 40.28 \\
\hline
\end{tabular}
\label{tab:errNOdiff}
\end{table}

\textit{Accuracy of the prediction.} Table~\ref{tab:errNOdiff} shows the errors of the $L^2$ projection onto each space and the errors of the prediction. From the constant-mode baseline $V_{ms,1}$ to the local POD space $V_{ms,20}$, the mean test error of the prediction decreases from 5.77\% to 1.20\% in the $L^2$ norm in Test~1, and from 5.69\% to 0.58\% in Test~2. For the global spaces $V_{pod,50}$ to $V_{pod,400}$ it decreases from 4.25\% to 1.31\% and from 1.41\% to 0.63\% in the $L^2$ norm. The mean training and test $L^2$ errors differ by at most 0.3 percentage points.
The network is trained towards the coefficients of the $L^2$ projection. For each realization its squared $L^2$ error is the sum of the squared projection and learning errors in \eqref{eq:split}. For the small spaces the mean prediction error is close to the mean projection error. In Test~1 their difference is 0.10 percentage points for $V_{ms,1}$ and 0.55 for $V_{ms,5}$; the corresponding energy entries are 68.8\% and 69.0\%, and 22.9\% and 21.8\%. For the larger spaces the difference between the mean $L^2$ errors grows to 0.93 and 1.08 percentage points for $V_{ms,10}$ and $V_{ms,20}$, while the projection error decreases further. In Test~2 the difference grows from 0.03 to 0.41 percentage points, and for the global spaces from 0.14 to 0.84 in Test~1. 
The limited improvement of the predictions in the larger spaces points to a bottleneck on the learning side.

\begin{table}[h!]
\centering
\footnotesize
\caption{Sparsity and memory of the reduced spaces on the $240\times240$ grid. The number of coarse degrees of freedom, the number of nonzero entries, the percentage of nonzero entries and the memory of the basis matrix $P$, and the same for the coarse matrix $A_H=P^TA(p)P$. }
\setlength{\tabcolsep}{4pt}
\begin{tabular}{lrrrrrrr}
\hline
 & & \multicolumn{3}{c}{Basis matrix $P$} & \multicolumn{3}{c}{Coarse matrix $A_H$} \\
 & DOF$_H$ & nonzeros & \% & MB & nonzeros & \% & MB \\
\hline
$V_{ms,1}$ & 81 & 230\,400 & 4.94 & 2.8 & 1\,325 & 20.2 & 0.02 \\
$V_{ms,5}$ & 405 & 1\,152\,000 & 4.94 & 13.8 & 33\,125 & 20.2 & 0.40 \\
$V_{ms,10}$ & 810 & 2\,304\,000 & 4.94 & 27.7 & 132\,500 & 20.2 & 1.59 \\
$V_{ms,20}$ & 1620 & 4\,608\,000 & 4.94 & 55.3 & 530\,000 & 20.2 & 6.37 \\
\hline
$V_{pod,50}$ & 50 & 2\,880\,000 & 100 & 23.0 & 2\,500 & 100 & 0.02 \\
$V_{pod,100}$ & 100 & 5\,760\,000 & 100 & 46.1 & 10\,000 & 100 & 0.08 \\
$V_{pod,200}$ & 200 & 11\,520\,000 & 100 & 92.2 & 40\,000 & 100 & 0.32 \\
$V_{pod,400}$ & 400 & 23\,040\,000 & 100 & 184.3 & 160\,000 & 100 & 1.28 \\
\hline
\end{tabular}
\label{tab:sparsity}
\end{table}

\textit{Local and global spaces.} With five modes per neighborhood, $V_{ms,5}$ gives smaller mean $L^2$ prediction errors than $V_{pod,50}$ and $V_{pod,100}$ in both tests. The space $V_{pod,400}$, with nearly the same number of coefficients, gives slightly smaller errors (1.31\% against 1.57\% and 0.63\% against 0.80\%), but its basis requires 13 times more memory. Table~\ref{tab:sparsity} compares the storage of the spaces. A local basis function is supported on one coarse neighborhood, so that the matrix $P$ of $V_{ms,m}$ has 4.94\% of nonzero entries, while the matrix $\Phi_r$ of a global space is dense. The basis of $V_{ms,5}$ requires 13.8 MB and the basis of $V_{pod,400}$, with nearly the same number of degrees of freedom, 184.3 MB (221 MB and 2.9 GB on the $960\times960$ grid). The coarse matrix of $V_{ms,m}$ has at most $21m$ nonzero entries per row. The nearly equal-dimensional global space also has the smaller $L^2$ projection error (0.47\% against 1.03\% in Test~1, Table~\ref{tab:errNOdiff}). The local spaces instead provide favorable representation accuracy relative to basis storage in these tests: $V_{ms,10}$ (27.7 MB) and $V_{pod,50}$ (23.0 MB) have projection errors of 0.34\% and 4.12\% in Test~1, and $V_{ms,20}$ (55.3 MB) and $V_{pod,200}$ (92.2 MB) have errors of 0.12\% and 1.21\%. 
These comparisons concern the stored basis matrices, not the total memory of the learning and correction pipeline. The memory is given in megabytes for double precision, with the compressed sparse row format for the local spaces and dense matrices for the global spaces.

\begin{table}[h!]
\centering
\footnotesize
\caption{Errors on different output grids. All errors use the $960\times960$ solution restricted to the current grid as the reference. The rows show the fine-grid solve, the $L^2$ projection (P+), the Galerkin coarse solve (S+) and the neural prediction (NO+) for $V_{ms,1}$, $V_{ms,5}$ and $V_{pod,50}$. 
}
\setlength{\tabcolsep}{4pt}
\begin{tabular}{lccccc}
\hline
 & $60\times60$ & $120\times120$ & $240\times240$ & $480\times480$ & $960\times960$ \\
\hline
\multicolumn{6}{c}{Test 1a} \\
\hline
Fine solve & 0.87 / 1.60 & 0.30 / 0.89 & 0.10 / 0.45 & 0.03 / 0.19 & 0.00 / 0.00 \\
\hline
P+$V_{ms,1}$ & 5.73 / 46.47 & 5.75 / 54.01 & 5.77 / 62.66 & 5.78 / 71.77 & 5.78 / 82.08 \\
S+$V_{ms,1}$ & 11.44 / 38.28 & 12.61 / 40.84 & 13.62 / 42.93 & 14.16 / 44.01 & 14.44 / 44.57 \\
NO+$V_{ms,1}$ & 5.72 / 44.99 & 5.84 / 54.12 & 5.88 / 62.96 & 5.88 / 73.70 & 5.89 / 87.04 \\
\hline
P+$V_{ms,5}$ & 2.09 / 24.00 & 1.96 / 27.48 & 1.98 / 31.93 & 1.98 / 38.40 & 1.98 / 45.41 \\
S+$V_{ms,5}$ & 3.03 / 14.99 & 3.52 / 16.06 & 3.93 / 17.64 & 4.30 / 19.71 & 4.55 / 20.82 \\
NO+$V_{ms,5}$ & 2.37 / 23.08 & 2.44 / 27.24 & 2.47 / 31.27 & 2.47 / 35.09 & 2.47 / 50.28 \\
\hline
P+$V_{pod,50}$ & 5.62 / 40.83 & 5.65 / 41.49 & 5.66 / 41.79 & 5.67 / 43.88 & 5.67 / 54.82 \\
S+$V_{pod,50}$ & 10.04 / 32.49 & 10.69 / 33.30 & 10.99 / 33.83 & 11.85 / 36.29 & 13.43 / 48.88 \\
NO+$V_{pod,50}$ & 5.68 / 40.01 & 5.76 / 41.47 & 5.78 / 41.81 & 5.79 / 43.98 & 5.79 / 54.91 \\
\hline
\multicolumn{6}{c}{Test 1b} \\
\hline
Fine solve & 0.48 / 1.11 & 0.17 / 0.61 & 0.06 / 0.30 & 0.02 / 0.13 & 0.00 / 0.00 \\
\hline
P+$V_{ms,1}$ & 5.44 / 44.99 & 5.61 / 54.87 & 5.68 / 66.62 & 5.70 / 79.75 & 5.70 / 94.78 \\
S+$V_{ms,1}$ & 11.09 / 36.99 & 13.13 / 39.97 & 14.87 / 42.20 & 15.98 / 43.39 & 16.65 / 44.00 \\
NO+$V_{ms,1}$ & 5.42 / 44.15 & 5.67 / 55.56 & 5.74 / 67.52 & 5.76 / 82.39 & 5.76 / 99.13 \\
\hline
P+$V_{ms,5}$ & 0.91 / 13.98 & 0.85 / 16.11 & 0.87 / 19.31 & 0.96 / 27.03 & 0.98 / 36.06 \\
S+$V_{ms,5}$ & 1.20 / 9.89 & 1.54 / 10.95 & 1.84 / 12.48 & 2.42 / 15.02 & 2.80 / 16.31 \\
NO+$V_{ms,5}$ & 1.37 / 14.82 & 1.42 / 17.68 & 1.44 / 20.89 & 1.46 / 24.75 & 1.47 / 44.97 \\
\hline
P+$V_{pod,50}$ & 3.66 / 27.29 & 3.69 / 28.43 & 3.71 / 28.86 & 3.72 / 31.83 & 3.72 / 48.24 \\
S+$V_{pod,50}$ & 5.75 / 24.11 & 6.13 / 24.81 & 6.30 / 25.28 & 7.12 / 28.52 & 8.93 / 45.67 \\
NO+$V_{pod,50}$ & 3.74 / 28.29 & 3.79 / 29.51 & 3.81 / 29.94 & 3.82 / 32.87 & 3.82 / 48.89 \\
\hline
\multicolumn{6}{c}{Test 2a} \\
\hline
Fine solve & 0.58 / 4.26 & 0.25 / 2.59 & 0.10 / 1.45 & 0.03 / 0.64 & 0.00 / 0.00 \\
\hline
P+$V_{ms,1}$ & 4.88 / 70.51 & 4.96 / 87.03 & 5.00 / 106.79 & 5.01 / 128.49 & 5.02 / 153.05 \\
S+$V_{ms,1}$ & 15.71 / 57.17 & 17.84 / 63.30 & 19.67 / 67.64 & 20.80 / 70.04 & 21.45 / 71.33 \\
NO+$V_{ms,1}$ & 4.78 / 68.40 & 4.96 / 87.11 & 5.00 / 106.75 & 5.03 / 130.85 & 5.03 / 157.38 \\
\hline
P+$V_{ms,5}$ & 0.68 / 25.22 & 0.48 / 27.51 & 0.43 / 29.34 & 0.45 / 31.83 & 0.45 / 36.64 \\
S+$V_{ms,5}$ & 1.96 / 14.09 & 2.45 / 16.38 & 2.98 / 18.52 & 3.23 / 20.33 & 3.38 / 21.25 \\
NO+$V_{ms,5}$ & 0.42 / 22.52 & 0.45 / 25.96 & 0.45 / 28.26 & 0.48 / 33.94 & 0.49 / 43.73 \\
\hline
P+$V_{pod,50}$ & 0.86 / 23.93 & 0.74 / 25.28 & 0.72 / 25.99 & 0.73 / 32.10 & 0.74 / 42.63 \\
S+$V_{pod,50}$ & 3.30 / 16.78 & 3.71 / 18.47 & 3.88 / 19.33 & 4.09 / 27.05 & 4.29 / 39.12 \\
NO+$V_{pod,50}$ & 0.70 / 22.86 & 0.72 / 25.11 & 0.73 / 26.25 & 0.75 / 32.49 & 0.76 / 42.97 \\
\hline
\multicolumn{6}{c}{Test 2b} \\
\hline
Fine solve & 1.26 / 6.10 & 0.53 / 3.59 & 0.21 / 2.02 & 0.07 / 0.91 & 0.00 / 0.00 \\
\hline
P+$V_{ms,1}$ & 6.32 / 134.15 & 6.36 / 146.13 & 6.41 / 161.64 & 6.42 / 179.46 & 6.42 / 200.91 \\
S+$V_{ms,1}$ & 45.79 / 75.65 & 55.62 / 83.95 & 62.76 / 88.77 & 67.27 / 91.44 & 69.99 / 92.85 \\
NO+$V_{ms,1}$ & 6.15 / 128.49 & 6.37 / 145.41 & 6.43 / 161.54 & 6.46 / 180.62 & 6.46 / 203.55 \\
\hline
P+$V_{ms,5}$ & 1.42 / 65.93 & 0.92 / 75.07 & 0.81 / 81.21 & 0.80 / 84.52 & 0.81 / 89.10 \\
S+$V_{ms,5}$ & 10.71 / 34.57 & 14.60 / 42.22 & 17.66 / 47.27 & 19.77 / 50.74 & 21.44 / 53.07 \\
NO+$V_{ms,5}$ & 0.87 / 62.56 & 0.94 / 75.18 & 0.96 / 82.13 & 0.97 / 86.54 & 0.98 / 92.66 \\
\hline
P+$V_{pod,50}$ & 1.96 / 71.99 & 1.68 / 78.89 & 1.62 / 82.37 & 1.63 / 86.62 & 1.63 / 92.46 \\
S+$V_{pod,50}$ & 15.64 / 42.07 & 18.41 / 47.41 & 19.49 / 50.07 & 20.32 / 55.08 & 20.75 / 63.22 \\
NO+$V_{pod,50}$ & 1.63 / 68.55 & 1.67 / 77.54 & 1.68 / 81.91 & 1.70 / 86.58 & 1.70 / 92.58 \\
\hline
\end{tabular}
\label{tab:errres}
\end{table}

\textit{Comparison with the projection and the Galerkin coarse solve.} In Table~\ref{tab:errres}, we compare the neural prediction with two other approximations in the same reduced space. The first is the $L^2$ projection of the fine-grid solution, $\Pi p=PU_*$, which is the best approximation in the $L^2$ norm. The second is the Galerkin coarse solution, $p_G=PU^H$, obtained by the Picard iteration \eqref{eq:coarse} with the tolerance $10^{-6}$. Unlike the projection, it does not require the reference solution $p$. 
The network is trained once with the data of the $240\times240$ grid. On another output grid the same predicted coefficients are used, and the pressure is reconstructed with the basis transferred to that grid (Section~\ref{sec:nn}). 
The permeability remains piecewise constant on the same input grid throughout these tests, so refining the output grid does not introduce new input information.
Table~\ref{tab:errres} shows the errors on the grids from $60\times60$ to $960\times960$ for the four separate tests. On each grid the fine solve, the $L^2$ projection of its solution and the Galerkin solution are computed on that grid, with the local modes transferred to the grid as in Section~\ref{sec:nn}. The prediction is computed once, reconstructed with the basis of the $240\times240$ grid and transferred to the grid by cell averages or by bilinear interpolation. All errors are computed on the current grid, and the reference is the solution of the $960\times960$ grid restricted to this grid.
The fine-grid solutions approach the common reference as the grid is refined. On the $240\times240$ grid of the training data their errors are between 0.06\% and 0.21\% in the $L^2$ norm and between 0.3\% and 2.0\% in the state-dependent energy measure. The $L^2$ error of the prediction NO+$V_{ms,5}$ changes by at most 0.05 percentage points on the grids from $120\times120$ to $960\times960$. In all tests and on all grids, except for $V_{ms,5}$ on the $60\times60$ grid in Test~1b, the prediction has a smaller $L^2$ error than the Galerkin solution in the same space. The space $V_{ms,5}$ gives smaller $L^2$ errors than $V_{ms,1}$ and $V_{pod,50}$ for the projection and for the prediction. 
The energy entries of the projection and prediction grow under refinement (from 23.1\% on the $60\times60$ grid to 50.3\% on the $960\times960$ grid for NO+$V_{ms,5}$ in Test~1a). Transferring the basis does not recover pressure scales unresolved on the $240\times240$ training grid and the growth can also depend on the boundary treatment in the transfer. The Galerkin solution is recomputed on every grid, and its $L^2$ error also grows (from 10.71\% to 21.44\% in Test~2b). Thus, these tests establish transfer of the decoded prediction without retraining, rather than generalization to arbitrary input resolutions or recovery of new fine-scale information.

\begin{table}[h!]
\centering
\footnotesize
\caption{Computational time (in seconds) for the standard fine-grid solve, neural operator prediction (NO+$V_{ms,5}$ and NO+$V_{pod,50}$), and the corresponding reduced-order solves on $V_{ms,5}$ and $V_{pod,50}$ (S+$V_{ms,5}$ and S+$V_{pod,50}$). 
Mean $\pm$ standard deviation over five test realizations. n.c. means that the AMG-CG iteration did not converge for at least one of the realizations.}
\setlength{\tabcolsep}{2pt}
\begin{tabular}{lccccc}
\hline
 & $60\times60$ & $120\times120$ & $240\times240$ & $480\times480$ & $960\times960$ \\
\hline
\multicolumn{6}{c}{Test 1}\\
\hline
Fine solve & 0.0674 $\pm$ 0.0054 & 0.3118 $\pm$ 0.0050 & 1.7611 $\pm$ 0.0370 & 11.454 $\pm$ 0.268 & 75.775 $\pm$ 2.182 \\
Fine solve, AMG-CG & 0.1716 $\pm$ 0.0083 & 0.6223 $\pm$ 0.0365 & 2.2553 $\pm$ 0.2211 & 10.094 $\pm$ 1.179 & 44.495 $\pm$ 2.393 \\
NO+$V_{pod,50}$ & 0.0010 $\pm$ 0.0001 & 0.0011 $\pm$ 0.0000 & 0.0009 $\pm$ 0.0000 & 0.0051 $\pm$ 0.0000 & 0.0172 $\pm$ 0.0001 \\
NO+$V_{ms,5}$ & 0.0014 $\pm$ 0.0001 & 0.0015 $\pm$ 0.0001 & 0.0013 $\pm$ 0.0001 & 0.0058 $\pm$ 0.0007 & 0.0179 $\pm$ 0.0007 \\
S+$V_{pod,50}$ & 0.0173 $\pm$ 0.0019 & 0.0550 $\pm$ 0.0014 & 0.2283 $\pm$ 0.0093 & 1.0089 $\pm$ 0.0283 & 3.9761 $\pm$ 0.0735 \\
S+$V_{ms,5}$ & 0.0847 $\pm$ 0.0027 & 0.2655 $\pm$ 0.0071 & 0.8753 $\pm$ 0.0319 & 3.1991 $\pm$ 0.0846 & 11.964 $\pm$ 0.277 \\
\hline
\multicolumn{6}{c}{Test 2}\\
\hline
Fine solve & 0.0519 $\pm$ 0.0094 & 0.2553 $\pm$ 0.0422 & 1.4245 $\pm$ 0.2399 & 9.3079 $\pm$ 1.5347 & 61.213 $\pm$ 9.918 \\
Fine solve, AMG-CG & 0.1999 $\pm$ 0.0845 & n.c. & n.c. & n.c. & n.c. \\
NO+$V_{pod,50}$ & 0.0010 $\pm$ 0.0000 & 0.0012 $\pm$ 0.0001 & 0.0010 $\pm$ 0.0000 & 0.0051 $\pm$ 0.0001 & 0.0172 $\pm$ 0.0001 \\
NO+$V_{ms,5}$ & 0.0013 $\pm$ 0.0000 & 0.0016 $\pm$ 0.0001 & 0.0013 $\pm$ 0.0001 & 0.0054 $\pm$ 0.0000 & 0.0175 $\pm$ 0.0002 \\
S+$V_{pod,50}$ & 0.0143 $\pm$ 0.0041 & 0.0431 $\pm$ 0.0074 & 0.1736 $\pm$ 0.0321 & 0.7612 $\pm$ 0.1302 & 3.1571 $\pm$ 0.5440 \\
S+$V_{ms,5}$ & 0.0729 $\pm$ 0.0154 & 0.2210 $\pm$ 0.0342 & 0.7071 $\pm$ 0.1114 & 2.6105 $\pm$ 0.4143 & 9.8341 $\pm$ 1.5399 \\
\hline
\end{tabular}
\label{tab:time1}
\end{table}

\textit{Computational time.} The prediction requires no nonlinear iteration and no assembly of the fine-grid system, unlike the fine-grid and the Galerkin coarse solves, and this difference grows as the grid is refined. Table~\ref{tab:time1} compares their online times (mean $\pm$ standard deviation over test realizations). The fine and Galerkin solves use the Picard tolerance $10^{-6}$. The fine solve uses a sparse direct solver in every Picard iteration, and the row AMG-CG is the same Picard iteration with conjugate gradients preconditioned by algebraic multigrid (smoothed aggregation, relative residual $10^{-8}$, at most $100$ iterations per linear solve). In Test~1 every linear solve reaches the tolerance. In Test~2, from the $120\times120$ grid on, at least one realization per grid fails to reach that linear tolerance under these AMG-CG settings (n.c.). The time of the neural operator includes the computation of the local input features, the evaluation of the network on the CPU and the reconstruction of the field on the grid. Its prediction is an approximate solution, whereas the solver timings correspond to the stated convergence tolerances.
The prediction time is below 0.02~s on all grids. On the $960\times960$ grid the prediction takes 0.018~s in Test~1 and 0.017~s in Test~2, the fine solve 76~s and 61~s with the direct solver and 44~s with AMG-CG in Test~1, and the Galerkin solve on $V_{ms,5}$ 12.0~s and 9.8~s. In Test~1 the prediction is thus about 4200 times faster than the fine solve with the direct solver, about 2500 times faster than the fine solve with AMG-CG and about 670 times faster than the Galerkin solve on the same space. Most of the time of the Galerkin solve on $V_{ms,5}$ is spent on the projection of the fine-grid matrix onto the coarse space in every Picard iteration.

\begin{table}[h!]
\centering
\footnotesize
\caption{Offline cost (first five rows) and online cost for the same number of permeability fields (last three rows). The times of the fine-grid solves, of the Galerkin solves and of the predictions are estimated with the mean time of one solve on the $240\times240$ grid (Table~\ref{tab:time1}). 
}
\setlength{\tabcolsep}{2pt}
\begin{tabular}{lcc}
\hline
 & Test 1 & Test 2 \\
\hline
Fine-grid solves & $\approx$ 8.3~h (17000 solves, $\approx$1.76~s each) & $\approx$ 4.7~h (12000 solves, $\approx$1.42~s each) \\
Local modes of $V_{ms,m}$ & 390~s & 589~s \\
Global modes of $V_{pod,r}$ & 33~s & 58~s \\
Training of NO+$V_{ms,5}$ & 418~s & 420~s \\
Training of NO+$V_{pod,50}$ & 403~s & 393~s \\
\hline
Galerkin solves on $V_{ms,5}$ & $\approx$ 4.1~h (17000 solves, $\approx$0.88~s each) & $\approx$ 2.4~h (12000 solves, $\approx$0.71~s each) \\
Predictions NO+$V_{ms,5}$ & $\approx$ 22~s (17000 predictions, $\approx$0.0013~s each) & $\approx$ 16~s (12000 predictions, $\approx$0.0013~s each) \\
Predictions NO+$V_{pod,50}$ & $\approx$ 16~s (17000 predictions, $\approx$0.0009~s each) & $\approx$ 12~s (12000 predictions, $\approx$0.0010~s each) \\
\hline
\end{tabular}
\label{tab:offline}
\end{table}

\textit{Offline and online cost.} Table~\ref{tab:offline} shows the offline and the online cost. It is dominated by the fine-grid solves of the datasets with the sparse direct solver (about 8.3~h and 4.7~h), while the local modes take 390~s and 589~s and the training of a network about 410~s for both spaces. The training is fast because it involves only coarse quantities. The network maps the $810$ input features ($q=10$ principal components of the log-permeability on each of the $81$ neighborhoods) to the coefficients of the reduced space ($405$ for $V_{ms,5}$), and the $57600$ unknowns of the fine grid do not appear in the training. 
The loss \eqref{eq:loss} is identical to the loss evaluated on the fine grid, and one training step is about five times faster in this implementation. 
The local modes of the four shape groups are computed independently of each other, and the local snapshots are extracted independently on the $81$ neighborhoods, so that this step can be performed in parallel. The projection error can be evaluated before network training when the reference solutions are available (column Projection of Table~\ref{tab:errNOdiff}). The remaining learning error depends on the input representation as well as on the network, its training and the data (Section~\ref{sec:erran}).
The last three rows show the estimated online cost for the same number of permeability fields. The Galerkin coarse solve requires, in every Picard iteration, the assembly of the fine-grid matrix at the current approximation, its projection onto the reduced space and the solution of the coarse system. On $V_{ms,5}$ this takes about 0.88~s and 0.71~s per realization, and about 4.1~h and 2.4~h for all realizations. The prediction requires only the input features, one evaluation of the network and one product with $P$, about 0.0013~s per realization, so that the corresponding predictions for $17000$ and $12000$ realizations take about 22~s and 16~s with $V_{ms,5}$ and 16~s and 12~s with $V_{pod,50}$. The prediction is thus more than 500 times faster online than the Galerkin coarse solve on the same space. 

\subsection{Correction of the Neural Prediction}
\label{sec:rescorr}

\textit{Setting.} The prediction has a small $L^2$ error, but its energy error can remain large (Table~\ref{tab:errNOdiff}). We therefore correct it on the fine grid in two ways, by the local correction \eqref{eq:localcorr}, with the coarse step in the space of the prediction, and by the global correction \eqref{eq:lincorr}, whose steps are solved by the conjugate gradient method with one V-cycle of smoothed-aggregation algebraic multigrid (AMG) as the preconditioner. The CG iteration of a step is stopped when the relative residual of its linear system is below $\varepsilon=10^{-2}$. Both corrections are started from the neural predictions on $V_{ms,1}$, $V_{ms,5}$ and $V_{pod,50}$. The complete fine-grid solve is given for comparison. It starts from the Darcy solution. 
All runs use the four fixed tests, the $240\times240$, $480\times480$ and $960\times960$ grids and four threads. The errors are computed on the given grid as in Table~\ref{tab:errres}, and the correction times contain the assembly, the setup of the preconditioner or of the local problems, and the iterations. 

\begin{figure}[h!]
\centering
\includegraphics[width=\linewidth]{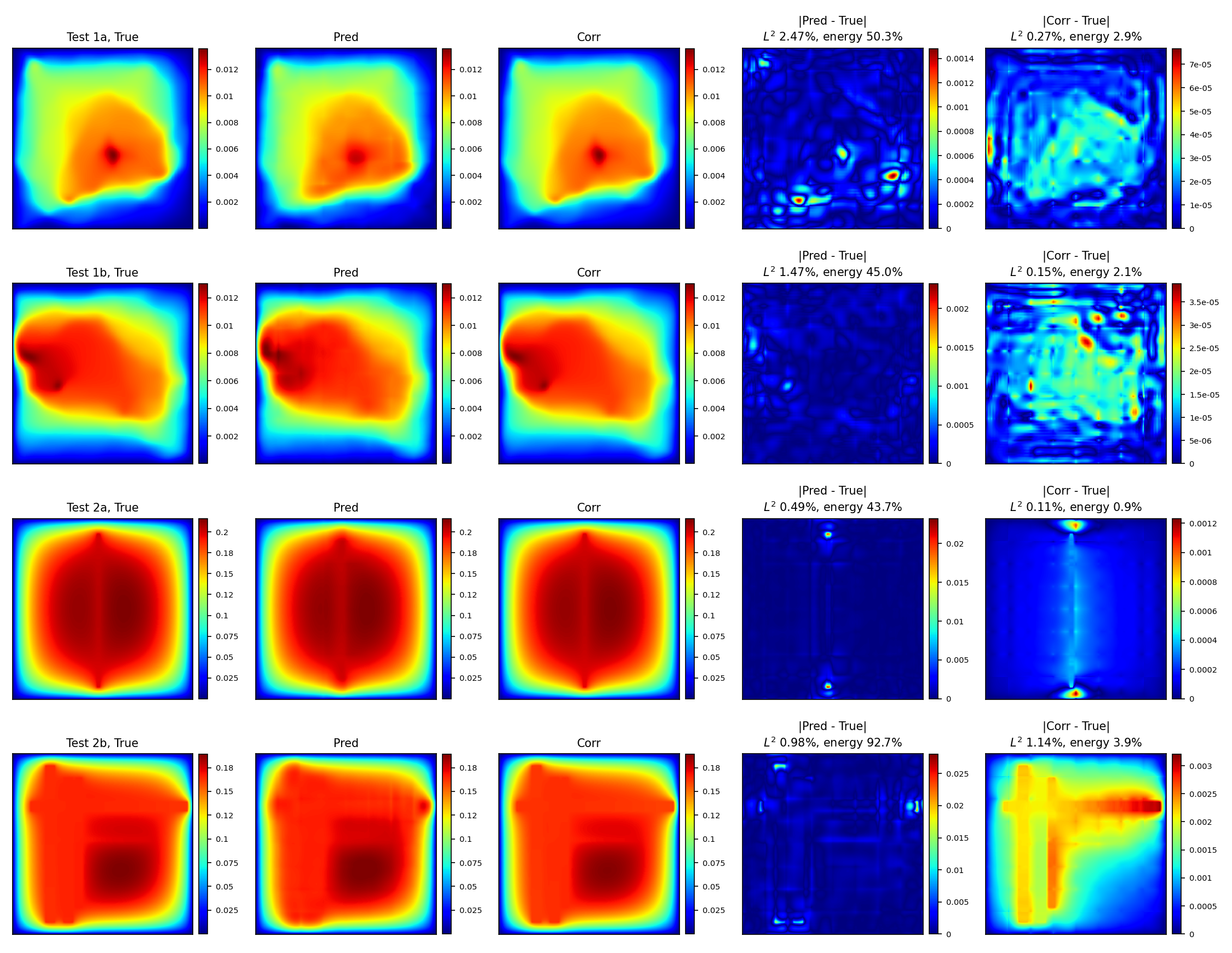}
\caption{Local correction on the $960\times960$ grid. The rows show Tests 1a, 1b, 2a and 2b, and the columns show the reference solution (True), the neural prediction NO+$V_{ms,5}$ (Pred), the solution after three steps of the local correction \eqref{eq:localcorr} (Corr) and the pointwise errors $|\mathrm{Pred}-\mathrm{True}|$ and $|\mathrm{Corr}-\mathrm{True}|$, with their $L^2$ and energy errors.}
\label{fig:corr960new}
\end{figure}

\textit{Errors of the local correction.} Figure~\ref{fig:corr960new} shows the solutions and the pointwise errors on the $960\times960$ grid for NO+$V_{ms,5}$ before and after three steps of the local correction. For the three predictions and the three grids, three steps reduce the energy error from between 43\% and 93\% to between 0.9\% and 4\%. The $L^2$ error after three steps is 0.1 to 1.14\% from NO+$V_{ms,5}$.

\begin{table}[h!]
\centering
\scriptsize
\caption{Offline time-to-accuracy evaluation of the local correction \eqref{eq:localcorr} and the global correction \eqref{eq:lincorr} with AMG-preconditioned CG. The correction rows give the first step at which the reference-based $L^2$ error is below $0.5\%$ and the state-dependent energy error is below $5\%$, and the associated wall time in seconds. 
The fine-solve rows report complete solves at Picard tolerance $10^{-6}$. Total CG iterations are given in parentheses. The entry n.r. means not reached in ten correction steps, and n.c. means failure to meet the linear tolerance under the stated AMG-CG settings.}
\label{tab:lincorr}
\setlength{\tabcolsep}{2.5pt}
\resizebox{\linewidth}{!}{\begin{tabular}{l|ccc|ccc|ccc|ccc}
\hline
 & \multicolumn{3}{c|}{Test 1a} & \multicolumn{3}{c|}{Test 1b} & \multicolumn{3}{c|}{Test 2a} & \multicolumn{3}{c}{Test 2b} \\
Start & $240^2$ & $480^2$ & $960^2$ & $240^2$ & $480^2$ & $960^2$ & $240^2$ & $480^2$ & $960^2$ & $240^2$ & $480^2$ & $960^2$ \\
\hline
\multicolumn{13}{c}{Number of steps} \\
\hline
\multicolumn{13}{l}{\textit{Local correction}} \\
NO+$V_{pod,50}$ & 6 & 6 & 7 & 4 & 5 & 5 & 2 & 2 & 2 & n.r. & n.r. & n.r. \\
NO+$V_{ms,1}$ & 7 & 7 & 8 & 5 & 6 & 6 & 5 & 6 & 7 & n.r. & n.r. & n.r. \\
NO+$V_{ms,5}$ & 2 & 2 & 3 & 2 & 2 & 2 & 2 & 2 & 2 & 3 & 4 & 6 \\
\multicolumn{13}{l}{\textit{Global correction with CG
}} \\
NO+$V_{pod,50}$ & 3 (12) & 3 (13) & 4 (20) & 2 (8) & 3 (14) & 3 (14) & 2 (14) & 2 (17) & 2 (23) & 2 (19) & 2 (19) & 2 (22) \\
NO+$V_{ms,1}$ & 3 (12) & 4 (18) & 4 (23) & 3 (12) & 3 (14) & 4 (20) & 4 (59) & 4 (63) & 4 (69) & 3 (36) & 4 (70) & 4 (80) \\
NO+$V_{ms,5}$ & 2 (7) & 2 (7) & 3 (13) & 2 (7) & 2 (7) & 3 (12) & 2 (10) & 2 (14) & 2 (12) & 2 (16) & 2 (21) & 2 (21) \\
\multicolumn{13}{l}{\textit{Fine solve, Picard iterations}} \\
direct & 16 & 16 & 16 & 16 & 16 & 16 & 15 & 15 & 15 & 14 & 14 & 14 \\
AMG-CG & 16 (249) & 16 (285) & 16 (361) & 16 (233) & 16 (259) & 16 (340) & n.c. & n.c. & n.c. & n.c. & n.c. & n.c. \\
\hline
\multicolumn{13}{c}{Wall time (in seconds)} \\
\hline
\multicolumn{13}{l}{\textit{Local correction}} \\
NO+$V_{pod,50}$ & 4.33 & 16.4 & 161.5 & 2.91 & 13.9 & 114.7 & 1.43 & 5.60 & 47.6 & n.r. & n.r. & n.r. \\
NO+$V_{ms,1}$ & 4.21 & 18.3 & 186.0 & 2.94 & 15.6 & 135.4 & 2.98 & 15.8 & 160.5 & n.r. & n.r. & n.r. \\
NO+$V_{ms,5}$ & 1.33 & 5.55 & 69.2 & 1.34 & 5.61 & 46.9 & 1.35 & 5.59 & 47.2 & 1.95 & 11.1 & 137.0 \\
\multicolumn{13}{l}{\textit{Global correction with CG}} \\
NO+$V_{pod,50}$ & 0.22 & 0.86 & 4.91 & 0.16 & 0.88 & 3.65 & 0.18 & 0.69 & 3.45 & 0.19 & 0.77 & 3.46 \\
NO+$V_{ms,1}$ & 0.22 & 1.15 & 5.10 & 0.24 & 0.91 & 5.05 & 0.54 & 2.00 & 8.63 & 0.35 & 2.17 & 9.44 \\
NO+$V_{ms,5}$ & 0.14 & 0.57 & 3.49 & 0.16 & 0.53 & 3.31 & 0.15 & 0.63 & 2.57 & 0.17 & 0.74 & 3.29 \\
\multicolumn{13}{l}{\textit{Fine solve}} \\
direct & 1.72 & 11.0 & 72.3 & 1.70 & 11.0 & 72.4 & 1.66 & 9.84 & 65.1 & 1.49 & 9.22 & 60.8 \\
AMG-CG & 2.09 & 8.86 & 42.5 & 2.14 & 8.66 & 41.9 & n.c. & n.c. & n.c. & n.c. & n.c. & n.c. \\
\hline
\end{tabular}}
\end{table}

\textit{Steps and time to the reference-based targets.} Table~\ref{tab:lincorr} shows the first correction step at which an $L^2$ error below $0.5\%$ and a state-dependent energy error below $5\%$ are reached, and the associated wall time. 
This is an offline evaluation using the known reference solution. 
Started from NO+$V_{ms,5}$, the global correction reaches these targets in 2 to 3 steps, with 7 to 13 CG iterations in total in Test~1 and 10 to 21 in Test~2. The complete fine-grid solver takes 14 to 16 Picard iterations to satisfy its convergence criterion. The correction targets are reached in 0.14 to 0.17~s on the $240\times240$ grid and in 2.6 to 3.5~s on the $960\times960$ grid, where the complete fine-grid solve with the direct solver takes 18 to 25 times longer. These timing ratios compare different accuracy requirements and reflect both neural initialization and the looser inner and outer accuracy requirements of the correction. The predictions on $V_{ms,1}$ and $V_{pod,50}$ have larger $L^2$ errors and need 3 to 4 and 2 to 4 steps, that is, 1.5 to 3.5 and 1.0 to 1.7 times the time from NO+$V_{ms,5}$.

The local correction reaches the targets in 2 to 6 steps from NO+$V_{ms,5}$, but every step solves all local problems, so that it takes 47 to 137~s on the $960\times960$ grid, approximately 14 to 42 times the global-correction time (8 to 42 times over all three grids). From NO+$V_{ms,1}$ and NO+$V_{pod,50}$ it reaches the targets within larger number of steps. Each step then needs more CG iterations and more time in all cases. In these experiments the neural prediction supplies an inexpensive starting state, and two or three global fine-grid correction steps reach the stated reference-based targets. Together with its sparse basis, this supports its use as an initialization with low basis-storage cost.

\section{Conclusions}
\label{sec:conclusions}

We developed a surrogate for nonlinear flow in random heterogeneous porous media. A neural operator predicts the coefficients of the solution in a local multiscale space that is built once from fine-grid snapshots. The prediction requires no nonlinear coarse solve and no approximation of the nonlinear terms, and the same coefficients can be decoded on other grids by transferring the basis.
On $V_{ms,5}$, the mean relative $L^2$ test errors are $1.57\%$ and $0.80\%$ for the KLE-based and the channelized media. They are close to the projection errors of the space ($1.03\%$ and $0.59\%$) and below the errors of the Galerkin solution on the same space.
The prediction takes less than $0.02$~s on all grids. On the $960\times960$ grid it is about $4200$ and $3500$ times faster than the fine-grid solve with the direct solver in Test~1 and Test~2, and about $670$ and $560$ times faster than the Galerkin coarse solve on the same space.
In the four representative correction tests, two or three AMG-preconditioned correction steps started from the prediction reach an $L^2$ error below $0.5\%$ and an energy error below $5\%$ in $2.6$ to $3.5$~s on the $960\times960$ grid. This is $18$ to $25$ times faster than the complete fine-grid solve with the direct solver.
The error of the prediction splits into a projection error and a learning error, and correction step gives computable indicators of both. For the larger spaces the learning error is the larger part, so that a further improvement requires more informative input features, a better network or more training data.

\bibliographystyle{plain}
\bibliography{lit}

@book{benner2017model,
  title={Model reduction and approximation: theory and algorithms},
  author={Benner, Peter and Ohlberger, Mario and Cohen, Albert and Willcox, Karen},
  year={2017},
  publisher={SIAM}
}

@book{quarteroni2015reduced,
  title={Reduced basis methods for partial differential equations: an introduction},
  author={Quarteroni, Alfio and Manzoni, Andrea and Negri, Federico},
  year={2015},
  publisher={Springer}
}

@article{chatterjee2000introduction,
  title={An introduction to the proper orthogonal decomposition},
  author={Chatterjee, Anindya},
  journal={Current science},
  pages={808--817},
  year={2000},
  publisher={JSTOR}
}

@book{hesthaven2016certified,
  title={Certified reduced basis methods for parametrized partial differential equations},
  author={Hesthaven, Jan S and Rozza, Gianluigi and Stamm, Benjamin},
  year={2016},
  publisher={Springer}
}

@article{quarteroni2011certified,
  title={Certified reduced basis approximation for parametrized partial differential equations and applications},
  author={Quarteroni, Alfio and Rozza, Gianluigi and Manzoni, Andrea},
  journal={Journal of Mathematics in Industry},
  volume={1},
  number={1},
  pages={3},
  year={2011},
  publisher={Springer}
}

@inproceedings{fanaskov2023spectral,
  title={Spectral neural operators},
  author={Fanaskov, Vladimir Sergeevich and Oseledets, Ivan V},
  booktitle={Doklady Mathematics},
  volume={108},
  number={Suppl 2},
  pages={S226--S232},
  year={2023},
  organization={Springer}
}

@article{rudikov2025locally,
  title={Locally Subspace-Informed Neural Operators for Efficient Multiscale PDE Solving},
  author={Rudikov, Alexander and Fanaskov, Vladimir and Stepanov, Sergei and Shan, Buzheng and Muravleva, Ekaterina and Efendiev, Yalchin and Oseledets, Ivan},
  journal={arXiv preprint arXiv:2505.16030},
  year={2025}
}

@article{oommen2022learning,
  title={Learning two-phase microstructure evolution using neural operators and autoencoder architectures},
  author={Oommen, Vivek and Shukla, Khemraj and Goswami, Somdatta and Dingreville, R{\'e}mi and Karniadakis, George Em},
  journal={npj Computational Materials},
  volume={8},
  number={1},
  pages={190},
  year={2022},
  publisher={Nature Publishing Group UK London}
}

@article{lu2021learning,
  title={Learning nonlinear operators via DeepONet based on the universal approximation theorem of operators},
  author={Lu, Lu and Jin, Pengzhan and Pang, Guofei and Zhang, Zhongqiang and Karniadakis, George Em},
  journal={Nature machine intelligence},
  volume={3},
  number={3},
  pages={218--229},
  year={2021},
  publisher={Nature Publishing Group UK London}
}

@article{li2020fourier,
  title={Fourier neural operator for parametric partial differential equations},
  author={Li, Zongyi and Kovachki, Nikola and Azizzadenesheli, Kamyar and Liu, Burigede and Bhattacharya, Kaushik and Stuart, Andrew and Anandkumar, Anima},
  journal={arXiv preprint arXiv:2010.08895},
  year={2020}
}

@article{kovachki2024operator,
  title={Operator learning: Algorithms and analysis},
  author={Kovachki, Nikola B and Lanthaler, Samuel and Stuart, Andrew M},
  journal={Handbook of Numerical Analysis},
  volume={25},
  pages={419--467},
  year={2024},
  publisher={Elsevier}
}

@article{zhang2025finite,
  title={Finite Element Representation Network (FERN) for Operator Learning with a Localized Trainable Basis},
  author={Zhang, Zecheng and Liu, Hao and Fu, Guosheng and Schaeffer, Hayden and Lin, Guang},
  journal={arXiv preprint arXiv:2510.26962},
  year={2025}
}

@article{lu2022comprehensive,
  title={A comprehensive and fair comparison of two neural operators (with practical extensions) based on fair data},
  author={Lu, Lu and Meng, Xuhui and Cai, Shengze and Mao, Zhiping and Goswami, Somdatta and Zhang, Zhongqiang and Karniadakis, George Em},
  journal={Computer Methods in Applied Mechanics and Engineering},
  volume={393},
  pages={114778},
  year={2022},
  publisher={Elsevier}
}

@article{wang2025reduced,
  title={Reduced-Basis Deep Operator Learning for Parametric PDEs with Independently Varying Boundary and Source Data},
  author={Wang, Yueqi and Lin, Guang},
  journal={arXiv preprint arXiv:2511.18260},
  year={2025}
}

@article{sharma2024ensemble,
  title={Ensemble and Mixture-of-Experts DeepONets for Operator Learning},
  author={Sharma, Ramansh and Shankar, Varun},
  journal={arXiv preprint arXiv:2405.11907},
  year={2024}
}

@article{mou2026neural,
  title={Neural-POD: A Plug-and-Play Neural Operator Framework for Infinite-Dimensional Functional Nonlinear Proper Orthogonal Decomposition},
  author={Mou, Changhong and Lu, Binghang and Lin, Guang},
  journal={arXiv preprint arXiv:2602.15632},
  year={2026}
}

@article{tripura2022wavelet,
  title={Wavelet neural operator: a neural operator for parametric partial differential equations},
  author={Tripura, Tapas and Chakraborty, Souvik},
  journal={arXiv preprint arXiv:2205.02191},
  year={2022}
}

@article{hou1997multiscale,
  title={A multiscale finite element method for elliptic problems in composite materials and porous media},
  author={Hou, Thomas Y and Wu, Xiao-Hui},
  journal={Journal of Computational Physics},
  volume={134},
  number={1},
  pages={169--189},
  year={1997},
  publisher={Elsevier}
}

@book{efendiev2009multiscale,
  title={Multiscale finite element methods: theory and applications},
  author={Efendiev, Yalchin and Hou, Thomas Y},
  volume={4},
  year={2009},
  publisher={Springer}
}

@article{efendiev2013generalized,
  title={Generalized multiscale finite element methods},
  author={Efendiev, Yalchin and Galvis, Juan and Hou, Thomas Y},
  journal={Journal of Computational Physics},
  volume={251},
  pages={116--135},
  year={2013},
  publisher={Elsevier}
}

@article{chung2014adaptive,
  title={An adaptive GMsFEM for high-contrast flow problems},
  author={Chung, Eric T and Efendiev, Yalchin and Li, Guanglian},
  journal={Journal of Computational Physics},
  volume={273},
  pages={54--76},
  year={2014},
  publisher={Elsevier}
}

@article{vasilyeva2020learning,
  title={Learning macroscopic parameters in nonlinear multiscale simulations using nonlocal multicontinua upscaling techniques},
  author={Vasilyeva, Maria and Leung, Wing T and Chung, Eric T and Efendiev, Yalchin and Wheeler, Mary},
  journal={Journal of Computational Physics},
  volume={412},
  pages={109323},
  year={2020},
  publisher={Elsevier}
}

@article{vasilyeva2021machine,
  title={Machine learning for accelerating macroscopic parameters prediction for poroelasticity problem in stochastic media},
  author={Vasilyeva, Maria and Tyrylgin, Aleksey},
  journal={Computers \& Mathematics with Applications},
  volume={84},
  pages={185--202},
  year={2021},
  publisher={Elsevier}
}

@article{wang2020deep,
  title={Deep multiscale model learning},
  author={Wang, Yating and Cheung, Siu Wun and Chung, Eric T and Efendiev, Yalchin and Wang, Min},
  journal={Journal of Computational Physics},
  volume={406},
  pages={109071},
  year={2020},
  publisher={Elsevier}
}

@article{vasilyeva2021preconditioning,
  title={Preconditioning Markov chain Monte Carlo method for geomechanical subsidence using multiscale method and machine learning technique},
  author={Vasilyeva, Maria and Tyrylgin, Aleksei and Brown, Donald L and Mondal, Anirban},
  journal={Journal of Computational and Applied Mathematics},
  volume={392},
  pages={113420},
  year={2021},
  publisher={Elsevier}
}

@article{ghasemi2016model,
  title={Model order reduction in porous media flow simulation using quadratic bilinear formulation},
  author={Ghasemi, Mohammadreza and Gildin, Eduardo},
  journal={Computational Geosciences},
  volume={20},
  number={3},
  pages={723--735},
  year={2016},
  publisher={Springer}
}

@inproceedings{chaturantabut2009discrete,
  title={Discrete empirical interpolation for nonlinear model reduction},
  author={Chaturantabut, Saifon and Sorensen, Danny C},
  booktitle={Proceedings of the 48h IEEE Conference on Decision and Control (CDC) held jointly with 2009 28th Chinese Control Conference},
  pages={4316--4321},
  year={2009},
  organization={IEEE}
}

@article{amsallem2012nonlinear,
  title={Nonlinear model order reduction based on local reduced-order bases},
  author={Amsallem, David and Zahr, Matthew J and Farhat, Charbel},
  journal={International Journal for Numerical Methods in Engineering},
  volume={92},
  number={10},
  pages={891--916},
  year={2012},
  publisher={Wiley Online Library}
}

@article{ghasemi2016localized,
  title={Localized model order reduction in porous media flow simulation},
  author={Ghasemi, Mohamadreza and Gildin, Eduardo},
  journal={Journal of Petroleum Science and Engineering},
  volume={145},
  pages={689--703},
  year={2016},
  publisher={Elsevier}
}

@article{vasilyeva2026implicit,
  title={An implicit-explicit scheme and adaptive multiscale approximation of generalized Forchheimer flow in fractured porous media},
  author={Vasilyeva, Maria and Mbroh, Nana Adjoah and Spiridonov, Denis and Timofeeva, Tatiana},
  journal={Computers \& Mathematics with Applications},
  volume={211},
  pages={1--25},
  year={2026},
  publisher={Elsevier}
}

@article{arraras2019geometric,
  title={Geometric multigrid methods for Darcy--Forchheimer flow in fractured porous media},
  author={Arrar{\'a}s, Andr{\'e}s and Gaspar, Francisco Jos{\'e} and Portero, Laura and Rodrigo, Carmen},
  journal={Computers \& Mathematics with Applications},
  volume={78},
  number={9},
  pages={3139--3151},
  year={2019},
  publisher={Elsevier}
}

@article{aulisa2009analysis,
  title={Analysis of generalized Forchheimer flows of compressible fluids in porous media},
  author={Aulisa, Eugenio and Bloshanskaya, Lidia and Hoang, Luan and Ibragimov, Akif},
  journal={Journal of Mathematical Physics},
  volume={50},
  number={10},
  year={2009},
  publisher={AIP Publishing}
}

@article{paszke2019pytorch,
  title={Pytorch: An imperative style, high-performance deep learning library},
  author={Paszke, Adam and Gross, Sam and Massa, Francisco and Lerer, Adam and Bradbury, James and Chanan, Gregory and Killeen, Trevor and Lin, Zeming and Gimelshein, Natalia and Antiga, Luca and others},
  journal={Advances in neural information processing systems},
  volume={32},
  year={2019}
}

@article{virtanen2020scipy,
  title={SciPy 1.0: fundamental algorithms for scientific computing in Python},
  author={Virtanen, Pauli and Gommers, Ralf and Oliphant, Travis E and Haberland, Matt and Reddy, Tyler and Cournapeau, David and Burovski, Evgeni and Peterson, Pearu and Weckesser, Warren and Bright, Jonathan and others},
  journal={Nature methods},
  volume={17},
  number={3},
  pages={261--272},
  year={2020},
  publisher={Nature Publishing Group US New York}
}

@article{lanthaler2022error,
  title={Error estimates for {DeepONets}: a deep learning framework in infinite dimensions},
  author={Lanthaler, Samuel and Mishra, Siddhartha and Karniadakis, George E},
  journal={Transactions of Mathematics and Its Applications},
  volume={6},
  number={1},
  pages={tnac001},
  year={2022}
}

@article{babuska1997partition,
  title={The partition of unity method},
  author={Babu{\v{s}}ka, Ivo and Melenk, Jens M},
  journal={International Journal for Numerical Methods in Engineering},
  volume={40},
  number={4},
  pages={727--758},
  year={1997}
}

@article{sirovich1987turbulence,
  title={Turbulence and the dynamics of coherent structures. {I}. {C}oherent structures},
  author={Sirovich, Lawrence},
  journal={Quarterly of Applied Mathematics},
  volume={45},
  number={3},
  pages={561--571},
  year={1987}
}

@book{ghanem1991stochastic,
  title={Stochastic Finite Elements: A Spectral Approach},
  author={Ghanem, Roger G and Spanos, Pol D},
  publisher={Springer},
  address={New York},
  year={1991}
}

@article{bhattacharya2021model,
  title={Model reduction and neural networks for parametric {PDE}s},
  author={Bhattacharya, Kaushik and Hosseini, Bamdad and Kovachki, Nikola B and Stuart, Andrew M},
  journal={The SMAI Journal of Computational Mathematics},
  volume={7},
  pages={121--157},
  year={2021}
}

@article{wang2019prediction,
  title={Prediction of discretization of {GMsFEM} using deep learning},
  author={Wang, Min and Cheung, Siu Wun and Chung, Eric T and Efendiev, Yalchin and Leung, Wing Tat and Wang, Yating},
  journal={Mathematics},
  volume={7},
  number={5},
  pages={412},
  year={2019},
  doi={10.3390/math7050412}
}

@article{wang2022deep,
  title={A deep learning based reduced order modeling for stochastic underground flow problems},
  author={Wang, Yiran and Chung, Eric and Fu, Shubin},
  journal={Journal of Computational Physics},
  volume={467},
  pages={111449},
  year={2022},
  doi={10.1016/j.jcp.2022.111449}
}

@article{kropfl2022operator,
  title={Operator compression with deep neural networks},
  author={Kr{\"o}pfl, Fabian and Maier, Roland and Peterseim, Daniel},
  journal={Advances in Continuous and Discrete Models},
  volume={2022},
  pages={29},
  year={2022},
  doi={10.1186/s13662-022-03702-y}
}

@article{kropfl2023neural,
  title={Neural network approximation of coarse-scale surrogates in numerical homogenization},
  author={Kr{\"o}pfl, Fabian and Maier, Roland and Peterseim, Daniel},
  journal={Multiscale Modeling \& Simulation},
  volume={21},
  number={4},
  pages={1457--1485},
  year={2023},
  doi={10.1137/22M1524278}
}

@article{liu2024learning,
  title={Learning a generalized multiscale prolongation operator},
  author={Liu, Yucheng and Fu, Shubin and Zhou, Yingjie and Ye, Changqing and Chung, Eric T},
  journal={arXiv preprint arXiv:2410.06832},
  year={2024}
}

@article{haines2026neural,
  title={A neural-network-based multiscale hybridizable discontinuous {G}alerkin method for solving {PDE}s in porous media},
  author={Haines, Tony and Shi, Ke},
  journal={arXiv preprint arXiv:2608.25850},
  year={2026}
}

@article{wang2022mosaic,
  title={Mosaic flows: A transferable deep learning framework for solving {PDE}s on unseen domains},
  author={Wang, Hengjie and Planas, Robert and Chandramowlishwaran, Aparna and Bostanabad, Ramin},
  journal={Computer Methods in Applied Mechanics and Engineering},
  volume={389},
  pages={114424},
  year={2022},
  doi={10.1016/j.cma.2021.114424}
}

@article{choi2025defining,
  title={Defining foundation models for computational science: A call for clarity and rigor},
  author={Choi, Youngsoo and Cheung, Siu Wun and Kim, Youngkyu and Tsai, Ping-Hsuan and Diaz, Alejandro N and Zanardi, Ivan and Chung, Seung Whan and Copeland, Dylan Matthew and Kendrick, Coleman and Anderson, William and Iliescu, Traian and Heinkenschloss, Matthias},
  journal={arXiv preprint arXiv:2505.22904},
  year={2025}
}

@article{chung2026latent,
  title={Latent space element method},
  author={Chung, Seung Whan and Choi, Youngsoo and Miller, Christopher and Springer, H Keo and Sullivan, Kyle T},
  journal={arXiv preprint arXiv:2601.01741},
  year={2026}
}

@article{kratsios2024mixtures,
  title={Mixtures of neural operators reduce active complexity in operator learning},
  author={Kratsios, Anastasis and Furuya, Takashi and Lara Benitez, Jose Antonio and Lassas, Matti and de Hoop, Maarten},
  journal={arXiv preprint arXiv:2404.09101},
  year={2024}
}

@inproceedings{stiller2020large,
  title={Large-scale neural solvers for partial differential equations},
  author={Stiller, Patrick and Bethke, Friedrich and B{\"o}hme, Maximilian and Pausch, Richard and Torge, Sunna and Debus, Alexander and Vorberger, Jan and Bussmann, Michael and Hoffmann, Nico},
  booktitle={Driving Scientific and Engineering Discoveries Through the Convergence of HPC, Big Data and AI (SMC 2020)},
  series={Communications in Computer and Information Science},
  pages={20--34},
  publisher={Springer},
  year={2020},
  doi={10.1007/978-3-030-63393-6_2}
}

@article{dehoop2022cost,
  title={The cost-accuracy trade-off in operator learning with neural networks},
  author={de Hoop, Maarten V and Huang, Daniel Zhengyu and Qian, Elizabeth and Stuart, Andrew M},
  journal={Journal of Machine Learning},
  volume={1},
  number={3},
  pages={299--341},
  year={2022},
  doi={10.4208/jml.220509}
}

@article{jacobs1991adaptive,
  title={Adaptive mixtures of local experts},
  author={Jacobs, Robert A and Jordan, Michael I and Nowlan, Steven J and Hinton, Geoffrey E},
  journal={Neural Computation},
  volume={3},
  number={1},
  pages={79--87},
  year={1991},
  doi={10.1162/neco.1991.3.1.79}
}

@inproceedings{shazeer2017outrageously,
  title={Outrageously large neural networks: The sparsely-gated mixture-of-experts layer},
  author={Shazeer, Noam and Mirhoseini, Azalia and Maziarz, Krzysztof and Davis, Andy and Le, Quoc and Hinton, Geoffrey and Dean, Jeff},
  booktitle={International Conference on Learning Representations},
  year={2017}
}

@book{bensoussan1978asymptotic,
  title={Asymptotic Analysis for Periodic Structures},
  author={Bensoussan, Alain and Lions, Jacques-Louis and Papanicolaou, George},
  publisher={North-Holland},
  address={Amsterdam},
  year={1978}
}

@article{durlofsky1991numerical,
  title={Numerical calculation of equivalent grid block permeability tensors for heterogeneous porous media},
  author={Durlofsky, Louis J},
  journal={Water Resources Research},
  volume={27},
  number={5},
  pages={699--708},
  year={1991},
  doi={10.1029/91WR00107}
}

@article{jenny2003multi,
  title={Multi-scale finite-volume method for elliptic problems in subsurface flow simulation},
  author={Jenny, Patrick and Lee, Seong H and Tchelepi, Hamdi A},
  journal={Journal of Computational Physics},
  volume={187},
  number={1},
  pages={47--67},
  year={2003},
  doi={10.1016/S0021-9991(03)00075-5}
}

@article{hughes1998variational,
  title={The variational multiscale method---a paradigm for computational mechanics},
  author={Hughes, Thomas J R and Feij{\'o}o, Gonzalo R and Mazzei, Luca and Quincy, Jean-Baptiste},
  journal={Computer Methods in Applied Mechanics and Engineering},
  volume={166},
  number={1-2},
  pages={3--24},
  year={1998},
  doi={10.1016/S0045-7825(98)00079-6}
}

@article{e2003heterogeneous,
  title={The heterogeneous multiscale methods},
  author={E, Weinan and Engquist, Bjorn},
  journal={Communications in Mathematical Sciences},
  volume={1},
  number={1},
  pages={87--132},
  year={2003},
  doi={10.4310/CMS.2003.v1.n1.a8}
}

@article{arbogast2007multiscale,
  title={A multiscale mortar mixed finite element method},
  author={Arbogast, Todd and Pencheva, Gergina and Wheeler, Mary F and Yotov, Ivan},
  journal={Multiscale Modeling \& Simulation},
  volume={6},
  number={1},
  pages={319--346},
  year={2007},
  doi={10.1137/060662587}
}

@article{malqvist2014localization,
  title={Localization of elliptic multiscale problems},
  author={M{\aa}lqvist, Axel and Peterseim, Daniel},
  journal={Mathematics of Computation},
  volume={83},
  number={290},
  pages={2583--2603},
  year={2014},
  doi={10.1090/S0025-5718-2014-02868-8}
}

@article{chung2018constraint,
  title={Constraint energy minimizing generalized multiscale finite element method},
  author={Chung, Eric T and Efendiev, Yalchin and Leung, Wing Tat},
  journal={Computer Methods in Applied Mechanics and Engineering},
  volume={339},
  pages={298--319},
  year={2018},
  doi={10.1016/j.cma.2018.04.010}
}

@article{hesthaven2018non,
  title={Non-intrusive reduced order modeling of nonlinear problems using neural networks},
  author={Hesthaven, Jan S and Ubbiali, Stefano},
  journal={Journal of Computational Physics},
  volume={363},
  pages={55--78},
  year={2018},
  doi={10.1016/j.jcp.2018.02.037}
}

\appendix

\section{Global Reduced Space}
\label{AppA}
Let $S$ be the centered snapshot matrix
\[
S=
\left[
p_1-\bar p,\;
p_2-\bar p,\;
\ldots,\;
p_{N_s}-\bar p
\right]
\in\mathbb{R}^{N\times N_s},
\]
where $\bar p=\frac{1}{N_s}\sum_{j}p_j$ is the mean of the snapshots.
The global POD basis is constructed by applying the singular value decomposition (SVD) to the snapshot matrix $S$ ($S S^T \varphi_k=\sigma_k^2\varphi_k$,  $\sigma_1\ge \sigma_2\ge \cdots \ge 0$).  

Keeping the first $r$ modes gives the reduced space
\[
V_{pod,r}
=
\operatorname{span}\{\varphi_1,\ldots,\varphi_r\},
\quad
\Phi_r = [\varphi_1,\ldots,\varphi_r]\in\mathbb{R}^{N\times r}.
\]

The errors of the $L^2$ projection and of the prediction in the spaces $V_{pod,r}$ are given in Table~\ref{tab:errNOdiff} together with the errors in the local spaces. Here the approximation is sought in the form $\bar p+\Phi_r c$.

The same global POD basis can also be used in an operator-learning model,
following the POD-DeepONet idea built on the DeepONet branch--trunk
architecture~\cite{lu2021learning,lu2022comprehensive}.  In this setting, the
trunk is fixed to the POD basis $\Phi_r$, and a neural network learns the
coefficient map, 
$\mathcal{G}_{\theta}^{pod}$ from the input features $x$ of Section~\ref{sec:nn} to the coefficients $\widehat c(x)\in\mathbb{R}^{r}$, and the prediction is $\widehat p=\bar p+\Phi_r\widehat c(x)$. 
The corresponding numerical results are presented in Section~\ref{sec:nolearn}.

\section{Error Analysis}
\label{AppB}

This appendix analyses the projection and the Galerkin coarse solution of Section~\ref{sec:coarse} and the neural prediction of Section~\ref{sec:no}. Throughout, we assume homogeneous Dirichlet boundary conditions, $a=0$, so that $b$ contains only the prescribed source contributions and is independent of the pressure. We also assume that the columns of $P$ are linearly independent; hence $M=P^TP$ is positive definite and the projectors below are well defined. The reduced space and the input encoder are fixed after their offline construction.

\subsection{Approximation Properties of the Multiscale Space}
\label{sec:approx}

The projection error is controlled by the eigenvalues that are discarded in the construction of the space.

\begin{lemma}
\label{prop:proj}
Let $\lambda^{g}_1\ge\lambda^{g}_2\ge\ldots$ be the eigenvalues of the matrix $\frac{1}{N_s}S_gS_g^T$ of the shape group $g$ in \eqref{eq:spod-pool} (if the modes are computed for each neighborhood separately, every group consists of one neighborhood and $S_g=S_I$). Then, for the snapshots from which the space is constructed,
\begin{equation}
\label{eq:tail}
\frac{1}{N_s}\sum_{j=1}^{N_s}\|(I-\Pi)p_j\|^2\le\sum_{g}\sum_{\ell>m}\lambda^{g}_{\ell} .
\end{equation}
\end{lemma}
\begin{proof}
Let $\Pi^m_g$ be the orthogonal projector onto the first $m$ modes of the group $g$, and let $w_j=\sum_g\sum_{I\in g}\chi^{\omega_I}R_I^T\Pi^m_gR_Ip_j\in V_{ms}$. Since $\chi^{\omega_I}\ge0$ and $\sum_I\chi^{\omega_I}=1$, the convexity of the square gives $\|p_j-w_j\|^2\le\sum_g\sum_{I\in g}\|(I-\Pi^m_g)R_Ip_j\|^2$, and $\frac{1}{N_s}\sum_j\sum_{I\in g}\|(I-\Pi^m_g)R_Ip_j\|^2=\sum_{\ell>m}\lambda^{g}_{\ell}$. Finally, $\|(I-\Pi)p_j\|\le\|p_j-w_j\|$.
\end{proof}

For the global POD space of dimension $r$ (Appendix~\ref{AppA}), \eqref{eq:tail} is replaced by the classical identity with the discarded eigenvalues of the global snapshot covariance. The bound \eqref{eq:tail} is a statement about the snapshots from which the space is constructed.

Under the homogeneous boundary condition, $A(v)=G^TT(v)G$ for every fine-grid vector $v$, where the matrix $G$ and the diagonal matrix $T(v)$ have one row for every face $e$ of the fine grid, with interior faces counted once. For the common face $e=E_{IJ}$ of the cells $K_I$ and $K_J$ we have $(Gv)_e=v_I-v_J$ and $T_e(v)=T_{IJ}(v)$, and for a face $e=E_{ID}$ of $K_I$ on $\partial\Omega$ we have $(Gv)_e=v_I$ and $T_e(v)=T_{ID}(v)$, with the positive transmissibilities of Section~\ref{sec:fine}. Thus $A(v)$ is positive definite. Let $A_G=A(p_G)$, let $\|v\|_{A_G}=(v^TA_Gv)^{1/2}$ be the energy norm of the problem with coefficients frozen at the Galerkin solution, and let $\Pi_G=P(P^TA_GP)^{-1}P^TA_G$ be the $A_G$-orthogonal projector onto $V_{ms}$.

\begin{proposition}
\label{prop:galerkin}
Let $p$ be the fine-grid solution, $A(p)\,p=b$, and let $p_G=PU^H$ be a solution of \eqref{eq:galerkin}. Then
\begin{equation}
\label{eq:galid}
\|p-p_G\|_{A_G}^2=\min_{v\in V_{ms}}\|p-v\|_{A_G}^2+\|\Pi_Gp-p_G\|_{A_G}^2 ,
\end{equation}
and
\begin{equation}
\label{eq:galnl}
\|\Pi_Gp-p_G\|_{A_G}\le\tau\,\|p\|_{A_G}\le\max_e\Big|1-\frac{T_e(p)}{T_e(p_G)}\Big|\,\|p\|_{A_G},\quad
\tau^2=\sum_e\omega_e\Big(1-\frac{T_e(p)}{T_e(p_G)}\Big)^2 ,
\end{equation}
where the sum is taken over the faces and $\omega_e=T_e(p_G)\,(Gp)_e^2/\|p\|_{A_G}^2\ge0$, $\sum_e\omega_e=1$, is the fraction of the energy $\|p\|_{A_G}^2$ on the face $e$. Thus $\tau$ is the root mean square of the relative change of the transmissibilities between $p_G$ and $p$, with the weights $\omega_e$, and it is bounded by the largest relative change. For a linear problem $\tau=0$, and $p_G=\Pi_Gp$ is the best approximation of $p$ in the energy norm.
\end{proposition}
\begin{proof}
We have $p-p_G=(I-\Pi_G)p+(\Pi_Gp-p_G)$, where the second term is in $V_{ms}$ and the first term is $A_G$-orthogonal to $V_{ms}$ and has the norm $\min_{v\in V_{ms}}\|p-v\|_{A_G}$. This gives \eqref{eq:galid}. By \eqref{eq:galerkin} and $A(p)\,p=b$, $P^TA_Gp_G=P^Tb=P^TA(p)\,p$, and therefore $P^TA_G(\Pi_Gp-p_G)=P^T(A_G-A(p))\,p$, that is, $\Pi_Gp-p_G=\Pi_Gz$ with $z=A_G^{-1}(A_G-A(p))\,p$. Since $\|\Pi_Gz\|_{A_G}\le\|z\|_{A_G}$ and, for every $v$,
\[
v^TA_Gz=\sum_e\big(T_e(p_G)-T_e(p)\big)(Gv)_e(Gp)_e\le\|v\|_{A_G}\Big(\sum_eT_e(p_G)\big(1-T_e(p)/T_e(p_G)\big)^2(Gp)_e^2\Big)^{1/2},
\]
we obtain \eqref{eq:galnl} with $\|p\|_{A_G}^2=\sum_eT_e(p_G)(Gp)_e^2$.
\end{proof}

Identity~\eqref{eq:galid} separates the best approximation error in the $A_G$ norm from a nonnegative nonlinear contribution. The latter is bounded by the energy-weighted relative mismatch of the transmissibilities in \eqref{eq:galnl}. If $\tau$ is small, the Galerkin solution is close to the energy projection, relative to $\|p\|_{A_G}$. 
In the $L^2$ norm, $\|(I-\Pi)p\|\le\|p-p_G\|$. When the transmissibility mismatch is small, $p_G$ is close to the energy projection $\Pi_Gp$, which can differ substantially from the $L^2$ projection $\Pi p$ in a heterogeneous medium. 
Moreover, an approximation trained towards the coefficients $U_*$ can therefore have a smaller $L^2$ error than the Galerkin solution in the same space, but cannot have a smaller $L^2$ error than $\Pi p$. 

\subsection{Error of the Neural Prediction}
\label{sec:erran}

In this section we give a decomposition of the unnormalised mean-square error of the neural prediction and correction-based indicators that can be computed without the reference solution. 
We use the notation of Section~\ref{sec:coarse}, 
the $L^2$-orthogonal projector $\Pi=PM^{-1}P^T$ onto the reduced space and the Euclidean norm $\|\cdot\|$ of fine-grid vectors. In addition, $\|V\|_M=(V^TMV)^{1/2}=\|PV\|$ for coefficient vectors. For a medium with the fine-grid solution $p$ and the input features $x$ of the network, the prediction is $\widehat p=P\mathcal{G}_{\theta}(x)$, and
\begin{equation}
\label{eq:split}
e=p-\widehat p,\quad e_{\Pi}=(I-\Pi)p,\quad e_L=\Pi(p-\widehat p)=P\big(U_*-\mathcal{G}_{\theta}(x)\big),\quad \|e\|^2=\|e_{\Pi}\|^2+\|e_L\|^2 ,
\end{equation}
where $e_{\Pi}$ is the projection error and $e_L$ is the learning error, $\|e_L\|=\|U_*-\mathcal{G}_{\theta}(x)\|_M$.

Let $\mathbb{E}$ denote expectation over new media drawn from the prescribed distribution, with the offline space and encoder held fixed. We assume that all second moments below are finite. Let $\bar p(x)=\mathbb{E}[\,p\,|\,x\,]$ be the conditional mean solution given the input features $x$.

\begin{theorem}
\label{prop:apriori}
For any network $\mathcal{G}_{\theta}$,
\begin{equation}
\label{eq:apriori}
\mathbb{E}\|p-P\mathcal{G}_{\theta}(x)\|^2
=\underbrace{\mathbb{E}\|(I-\Pi)p\|^2}_{\mathcal{E}_{\Pi}\text{ (projection)}}
+\underbrace{\mathbb{E}\|\Pi(p-\bar p(x))\|^2}_{\mathcal{E}_{x}\text{ (input compression)}}
+\underbrace{\mathbb{E}\|\Pi\bar p(x)-P\mathcal{G}_{\theta}(x)\|^2}_{\mathcal{E}_{\theta}\text{ (network)}} ,
\end{equation}
and the mean squared learning error is $\mathcal{E}_L=\mathbb{E}\|e_L\|^2=\mathcal{E}_{x}+\mathcal{E}_{\theta}$. Moreover, the following two statements hold. (a) $\mathcal{E}_{x}\le\mathbb{E}\|p-\bar p(x)\|^2$ for every reduced space. For fixed input features, $\mathcal{E}_{x}$ does not decrease along nested enlargements of the reduced space, and $\mathcal{E}_{x}=0$ if the features determine the medium. (b) For any map $x\mapsto\widehat q(x)$ into the fine-grid space, the coefficients $M^{-1}P^T\widehat q(x)$ have the mean squared learning error $\mathbb{E}\|\Pi(p-\widehat q(x))\|^2\le\mathbb{E}\|p-\widehat q(x)\|^2$.
\end{theorem}
\begin{proof}
The first term follows from \eqref{eq:split}. Further, $\Pi p-P\mathcal{G}_{\theta}(x)=\Pi(p-\bar p(x))+(\Pi\bar p(x)-P\mathcal{G}_{\theta}(x))$, where the second summand depends on $x$ only and $\mathbb{E}[\,\Pi(p-\bar p(x))\,|\,x\,]=0$, so that the mean of the scalar product of the two summands vanishes. Statements (a) and (b) follow from $\|\Pi v\|\le\|v\|$, from $\|\Pi_{V}v\|\le\|\Pi_{V'}v\|$ for $V\subset V'$, and from $p=\bar p(x)$ when $x$ determines the medium.
\end{proof}

The decomposition \eqref{eq:apriori} guides the choice of the reduced space. The projection error $\mathcal{E}_{\Pi}$ does not increase when modes are added to the space. Lemma~\ref{prop:proj} bounds it on the snapshots from which the space is constructed, but not on new media. The input-compression term $\mathcal{E}_{x}$ does not decrease when modes are added, and the network term $\mathcal{E}_{\theta}$ is not known a priori. It depends on the network, the training and the number of training samples, and error analyses for branch--trunk networks are given in \cite{lanthaler2022error,kovachki2024operator}. Statement (b) holds for any map, but it does not guarantee that a given network learns the coefficients as accurately as the fine-grid field. The number of modes is therefore chosen empirically, by balancing the projection error and the learning error. If the projection error is already small, a further improvement requires more informative input features, a better network or more training data. The decomposition is stated for the unnormalised squared error, while the loss \eqref{eq:loss} and the errors in the tables are relative.

In \eqref{eq:split} the learning error is measured in the norm $\|\cdot\|_M$ of the loss function. In the Euclidean norm of the coefficients one only has $\lambda_{\min}(M)\|V\|_2^2\le\|PV\|^2\le\lambda_{\max}(M)\|V\|_2^2$, so that a loss in the Euclidean norm would weight the error components by the eigenvalues of $M$. If every correction step reduces the error by a factor of at most $\bar\rho<1$, that is, $\|p-p^{k+1}\|\le\bar\rho\|p-p^k\|$ for every step and almost every medium, then the mean squared error after $n$ steps is bounded by $\bar\rho^{\,2n}(\mathcal{E}_{\Pi}+\mathcal{E}_L)$. This contraction is an assumption and is not proved here. The errors in \eqref{eq:split} require the reference solution. For a new medium, one correction step~\eqref{eq:lincorr} started from the prediction gives computable indicators. We first consider an exact linear solve and let
\begin{equation}
\label{eq:apostdef}
A(\widehat p)\,\delta=b-A(\widehat p)\,\widehat p,\quad \widehat p^{\,1}=\widehat p+\delta,\quad
\eta=\|\delta\|,\quad \eta_L=\|\Pi\delta\|=\|M^{-1}P^T\delta\|_M,\quad \eta_{\Pi}=\|(I-\Pi)\delta\|,
\end{equation}
so that $\eta^2=\eta_L^2+\eta_{\Pi}^2$. Here $\eta_L$ measures the part of the correction that lies in the reduced space, and $\eta_{\Pi}$ the part that the space cannot represent.

\begin{theorem}
\label{prop:apost}
Assume $\|p-\widehat p\|>0$, and let $\rho=\|p-\widehat p^{\,1}\|/\|p-\widehat p\|$ be the error reduction factor of the step. Then
\begin{equation}
\label{eq:apost}
(1-\rho)\,\|e\|\le\eta\le(1+\rho)\,\|e\|,
\quad
\big(\|e_L\|-\eta_L\big)^2+\big(\|e_{\Pi}\|-\eta_{\Pi}\big)^2\le\rho^2\|e\|^2 .
\end{equation}
In particular, if $\rho\le\bar\rho<1$, then
\[
\frac{\eta}{1+\bar\rho}\le\|e\|\le\frac{\eta}{1-\bar\rho},\quad
\big|\,\|e_L\|-\eta_L\big|\le\frac{\bar\rho}{1-\bar\rho}\,\eta,\quad
\big|\,\|e_{\Pi}\|-\eta_{\Pi}\big|\le\frac{\bar\rho}{1-\bar\rho}\,\eta,\quad
\|p-\widehat p^{\,1}\|\le\frac{\bar\rho}{1-\bar\rho}\,\eta .
\]
The first estimate in \eqref{eq:apost} holds in any norm.
\end{theorem}
\begin{proof}
Since $e=\delta+(p-\widehat p^{\,1})$, the triangle inequality gives the first estimate. Applying $\Pi$ and $I-\Pi$ and using $(I-\Pi)\widehat p=0$, we obtain $e_L=\Pi\delta+\Pi(p-\widehat p^{\,1})$ and $e_{\Pi}=(I-\Pi)\delta+(I-\Pi)(p-\widehat p^{\,1})$. Hence $|\,\|e_L\|-\eta_L|\le\|\Pi(p-\widehat p^{\,1})\|$ and $|\,\|e_{\Pi}\|-\eta_{\Pi}|\le\|(I-\Pi)(p-\widehat p^{\,1})\|$, and the squares of the two right-hand sides add up to $\|p-\widehat p^{\,1}\|^2=\rho^2\|e\|^2$.
\end{proof}

The indicators $\eta$, $\eta_L$ and $\eta_{\Pi}$ are computed from the prediction alone, without the reference solution and the training data, and the same construction applies to any approximation in the reduced space, including the Galerkin coarse solution. For a known bound $\bar\rho<1$ of the reduction factor, $\eta$ estimates the error, and $\eta_L$ and $\eta_{\Pi}$ estimate the learning error and the projection error, up to the deviations in Theorem~\ref{prop:apost}. The learning error is the sum of the input-compression error and the network error, and $\eta_L$ estimates the two together. In practice the reduction factor is estimated, at the cost of a second step, by the ratio of the norms of two successive corrections, $\widehat\rho=\|\delta^1\|/\|\delta^0\|$. This ratio is an empirical estimate of $\rho$ and is used in place of the bound $\bar\rho$. For an exact step, the argument of Proposition~\ref{prop:galerkin} also gives $\|p-\widehat p^{\,1}\|_{A(\widehat p)}\le\tau(\widehat p)\,\|p\|_{A(\widehat p)}$, where $\tau(\widehat p)$ is defined by \eqref{eq:galnl} with $p_G$ replaced by $\widehat p$ and measures the relative change of the transmissibilities between $\widehat p$ and $p$. The corrected solution leaves $V_{ms}$ and depends on the space only through the starting state $\widehat p$. The proof of \eqref{eq:apost} uses only the computed correction and its error reduction factor, so that the theorem also holds for an inexact linear solve, with $\rho$ and $\bar\rho$ of the computed step. When the prediction is the starting state of the fine-grid corrections, the first correction is available from the solver, and its splitting costs one additional application of $\Pi$.

\end{document}